# Observer design for nonlinear triangular systems with unobservable linearization

Dimitris Boskos[a] and John Tsinias[b]*

[a,b]Department of Mathematics, National Technical University of Athens, Zografou Campus 15780, Athens, Greece

The paper deals with the observer design problem for a wide class of triangular time-varying nonlinear systems, with unobservable linearization. Sufficient conditions are derived for the existence of a Luenberger-type observer, when it is a priori known that the initial state of the system belongs to a given nonempty bounded subset of the state space. For the general case, the state estimation is exhibited by means of a switching sequence of time-varying dynamics

## 1. Introduction

Observer design for nonlinear systems constitutes a problem of intense and growing research during the last two decades. Several approaches have been used for the solvability of this problem concerning a variety of systems, including those, whose dynamics have triangular structure (see e.g. Ahmed-Ali and Lamnabhi-Lagarrigue 1999; Ahrens and Khalil 2009; Alamir 1999; Andrieu and Praly 2006; Andrieu, Praly and Astolfi 2009; Arcak and Kokotovic 2001; Astolfi and Praly 2006; Back and Seo 2004, 2008; Besancon and Ticlea 2007; Gauthier, Hammouri and Othman 1992; Gauthier and Kupka 2001; Hammami 2005; Hammouri, Bornard and Busawon 2010; Jouan 2003; Karafyllis and Jiang 2011; Karafyllis and Kravaris 2009; Karagiannis, Carnevale and Astolfi, 2008; Kazantzis and Kravaris 1998; Kreisselmeier and Engel 2003; Krener and Xiao 2002a, 2002b; Liu 1997; Respondek, Pogromsky and Nijmeijer 2004; Sanfelice and Praly 2011; Starkov 1991; Tsinias 1990, 2008). Also, several approaches have been developed for the solvability of the related problem of output feedback stabilization by means of a dynamic controller; (see e.g. Krishnamurthy and Khorrami 2004; Qian and Lin 2006; Teel and Praly 1995; Tsinias 2000; Yang and Lin 2004).

In this work, we extent the methodology employed in Tsinias (2008), in order to show that the observer design problem is solvable for the following class of time-varying triangular systems:

$$\begin{aligned}
\dot{x}_1 &= f_1(t, x_1) + a_1(t, x_1) x_2^{m_1} \\
\dot{x}_2 &= f_2(t, x_1, x_2) + a_2(t, x_1) x_3^{m_2} \\
&\vdots \\
\dot{x}_{n-1} &= f_{n-1}(t, x_1, \ldots, x_{n-1}) + a_{n-1}(t, x_1) x_n^{m_{n-1}} \\
\dot{x}_n &= f_n(t, x_1, \ldots, x_n), \quad x_i \in \mathbb{R}, \quad i = 1, \ldots, n
\end{aligned} \quad (1a)$$

$$\text{with output } y = x_1, \quad (1b)$$

where it is assumed that

$$m_i, \ i = 1, \ldots, n-1 \text{ are odd positive integers} \quad (1c)$$

*Corresponding Author. Email: jtsin@central.ntua.gr

and under the additional assumptions that the values of the functions $|a_i(\cdot,\cdot)|$, $i=1,...,n-1$ are *strictly positive* and (1a) is *forward complete*. The paper is organized as follows. In the present section we first provide notations and various concepts, including the concept of the *switching observer,* that has been originally introduced in Boskos and Tsinias (2011), for general time-varying systems:

$$\dot{x} = f(t,x) \qquad (2a)$$

$$y = h(t,x) \qquad (2b)$$

$$(t,x) \in \mathbb{R}_{\geq 0} \times \mathbb{R}^n, \; y \in \mathbb{R}^{\bar{n}}$$

where $y(\cdot)$ plays the role of output. We then provide the precise statement of our main result (Theorem 1.1) concerning solvability of the observer design problem for (1), together with an illustrative example. Section 2 contains some preliminary results concerning solvability of the observer design problem for the case (2) with linear output (Propositions 2.1 and 2.2). In Section 3 we use the results of Section 2, in order to prove Theorem 1.1.

*Notations and definitions:* Throughout the paper we adopt the following notations. For a given vector $x \in \mathbb{R}^n$, $x'$ denotes its transpose and $|x|$ its Euclidean norm. We use the notation $|A| := \max\{|Ax| : x \in \mathbb{R}^n; |x|=1\}$ for the induced norm of a matrix $A \in \mathbb{R}^{m \times n}$. By $N$ we denote the class of all increasing $C^0$ functions $\phi : \mathbb{R}_{\geq 0} \to \mathbb{R}_{\geq 0}$. For given $R > 0$, we denote by $B_R$ the closed ball of radius $R > 0$, centered at $0 \in \mathbb{R}^n$. Consider a pair of metric spaces $X_1$, $X_2$ and a set-valued map $X_1 \ni x \to Q(x) \subset X_2$. We say that $Q(\cdot)$ satisfies the *Compactness Property* (**CP**), if for every sequence $(x_\nu)_{\nu \in \mathbb{N}} \subset X_1$ and $(q_\nu)_{\nu \in \mathbb{N}} \subset X_2$ with $x_\nu \to x \in X_1$ and $q_\nu \in Q(x_\nu)$, there exist a subsequence $(x_{\nu_k})_{k \in \mathbb{N}}$ and $q \in Q(x)$ such that $q_{\nu_k} \to q$. We also invoke the well known fact, see Karafyllis (2005), that the time-varying system (2a) is forward complete, if and only if there exists a function $\beta \in NN$ such that the solution $x(\cdot) := x(\cdot, t_0, x_0)$ of (2a) initiated from $x_0$ at time $t = t_0$ satisfies:

$$|x(t)| \leq \beta(t, |x_0|), \; \forall t \geq t_0 \geq 0, \; x_0 \in \mathbb{R}^n \qquad (3)$$

provided that the dynamics of (2a) are $C^0$ and (locally) Lipschitz on $x \in \mathbb{R}^n$. Under these regularity assumptions plus forward completeness for (2a), it follows that for each $t_0 \geq 0$ and $x_0 \in \mathbb{R}^n$ the corresponding output $y(t) = h(t, x(t, t_0, x_0))$ of (2) is defined for all $t \geq t_0$. Under previous assumptions, for each $t_0 \geq 0$ and nonempty subset $M$ of $\mathbb{R}^n$, we may consider the set $O(t_0, M)$ of all outputs of (2), corresponding to initial state $x_0 \in M$ and initial time $t_0 \geq 0$:

$$O(t_0, M) := \{ y : [t_0, \infty) \to \mathbb{R}^{\bar{n}} : y(t) = h(t, x(t, t_0, x_0)) \; ; \; x_0 \in M \} \qquad (4)$$

For given $\varnothing \neq M \subset \mathbb{R}^n$, we say that the **Observer Design Problem (ODP) is solvable for (2) with respect to** $M$, if for every $t_0 \geq 0$ there exist a continuous map $G := G_{t_0}(t, z, w) : [t_0, \infty) \times \mathbb{R}^n \times \mathbb{R}^{\bar{n}} \to \mathbb{R}^n$ and a nonempty set $\bar{M} \subset \mathbb{R}^n$ such that for every $z_0 \in \bar{M}$ and output $y \in O(t_0, M)$ the corresponding trajectory $z(\cdot) := z(\cdot, t_0, z_0)$; $z(t_0) = z_0$ of the observer $\dot{z} = G(t, z, y)$ exists for all $t \geq t_0$ and the error $e(t) := x(t) - z(t)$ between the trajectory $x(\cdot) := x(\cdot, t_0, x_0)$, $x_0 \in M$ of (2a) and the trajectory $z(\cdot) := z(\cdot, t_0, z_0)$ of the observer satisfies:

$$\lim_{t \to \infty} e(t) = 0 \tag{5}$$

We say that the **Switching Observer Design Problem (SODP) is solvable for (2) with respect to** $M$, if for every $t_0 \geq 0$ there exist:
- a strictly increasing sequence of times $(t_k)_{k \in \mathbb{N}}$ with

$$t_1 = t_0 \text{ and } \lim_{k \to \infty} t_k = \infty \tag{6}$$

- a sequence of continuous mappings

$$G_k := G_{k, t_{k-1}}(t, z, w) : [t_{k-1}, t_{k+1}] \times \mathbb{R}^n \times \mathbb{R}^{\bar{n}} \to \mathbb{R}^n, \ k \in \mathbb{N} \tag{7}$$

and a nonempty set $\bar{M} \subset \mathbb{R}^n$ such that the solution $z_k(\cdot)$ of the system

$$\dot{z}_k = G_k(t, z_k, y), \ t \in [t_{k-1}, t_{k+1}] \tag{8}$$

with initial $z(t_{k-1}) \in \bar{M}$ and output $y \in O(t_0, M)$, is defined for every $t \in [t_{k-1}, t_{k+1}]$ and in such a way that, if we consider the piecewise continuous map $Z : [t_0, \infty) \to \mathbb{R}^n$ defined as:

$$Z(t) := z_k(t), \ t \in [t_k, t_{k+1}), \ k \in \mathbb{N} \tag{9}$$

where for each $k \in \mathbb{N}$ the map $z_k(\cdot)$ denotes the solution of (8), then the error $e(t) := x(t) - Z(t)$ between the trajectory $x(\cdot) := x(\cdot, t_0, x_0)$, $x_0 \in M$ of (2a) and $Z(\cdot)$ satisfies (5).

Our main result is the following theorem.

**Theorem 1.1:** *Consider the system (1) and assume that $f_i \in C^1(\mathbb{R}_{\geq 0} \times \mathbb{R}^i; \mathbb{R})$, $i = 1, 2, ..., n$ and $a_i \in C^1(\mathbb{R}_{\geq 0} \times \mathbb{R}; \mathbb{R})$, $i = 1, 2, ..., n-1$. We assume that (1a) is forward complete; particularly, assume that the solution $x(\cdot) := x(\cdot, t_0, x_0)$ of (1a) satisfies the estimation (3) for certain $\beta \in NN$. Moreover, assume that*

$$|a_i(t, y)| > 0, \ \forall t \geq 0, \ y \in \mathbb{R}, \ i = 1, 2, ..., n-1 \tag{10}$$

*Then*
*(i) the SODP is solvable for (1) with respect to $\mathbb{R}^n$.*
*(ii) if in addition we assume that it is a priori known, that the initial states of (1) belong to a given ball $B_R$ of radius $R > 0$ centered at zero $0 \in \mathbb{R}^n$, then the ODP is solvable for (1) with respect to $B_R$.*

**Example 1.1:** Consider the 3-dimensional system $\dot{x}_1 = x_2^3 + f_1(t, x_1)$, $\dot{x}_2 = (1 + x_1^2)^{-1} x_3^3$ $+ f_2(t, x_1, x_2)$, $\dot{x}_3 = f_3(t, x_1, x_2, x_3)$ with output $y = x_1$, where $f_i \in C^1(\mathbb{R}_{\geq 0} \times \mathbb{R}^i; \mathbb{R})$, $i = 1, 2, 3$ and assume that there exists a function $C \in N$ such that for all $t \geq 0$, $x_1, x_2, x_3 \in \mathbb{R}$ it holds $x_1 f_1(t, x_1) \leq C(t) |x_1|^2$, $x_2 f_2(t, x_1, x_2) \leq C(t) |x_2|^2$ and $x_3 f_3(t, x_1, x_2, x_3) \leq C(t) |x_3|^2$. Clearly, the system has the form (1), with $m_1 = m_2 := 3$, and $a_1 := 1$, $a_2 := (1 + x_1^2)^{-1}$.

Also, by employing rest hypotheses, we can easily establish that it is forward complete. Therefore, all assumptions of Theorem 1.1 are fulfilled, hence, the SODP is solvable with respect to $\mathbb{R}^3$. If we further assume that it is a priori known, that the initial states belong to a given ball $B_R$ centered at zero $0 \in \mathbb{R}^3$, then the ODP is solvable with respect to $B_R$.

## 2. Preliminary results

In order to establish our main result concerning the case (1), we first need to prove some technical results concerning the case of systems (2) with linear output:

$$\dot{x} = f(t,x) \coloneqq F(t,x,H(t)x), \ (t,x) \in \mathbb{R}_{\geq 0} \times \mathbb{R}^n \quad (11a)$$

$$y = h(t,x) \coloneqq H(t)x, \ y \in \mathbb{R}^{\bar{n}} \quad (11b)$$

where $H : \mathbb{R}_{\geq 0} \to \mathbb{R}^{\bar{n} \times n}$ is $C^0$, $F : \mathbb{R}_{\geq 0} \times \mathbb{R}^n \times \mathbb{R}^{\bar{n}} \to \mathbb{R}^n$ is $C^0$ and (locally) Lipschitz on $(x,y) \in \mathbb{R}^n \times \mathbb{R}^{\bar{n}}$ and we assume that system (11a) is forward complete, namely, the solution $x(\cdot) \coloneqq x(\cdot, t_0, x_0)$ of (11a) satisfies (3) for certain $\beta \in NN$. We next adopt the following notation. For given $R > 0$ and $t \geq 0$ we define:

$$Y_R(t) \coloneqq \{ y \in \mathbb{R}^{\bar{n}} : y = H(t)x(t,t_0,x_0), \text{ for certain } t_0 \in [0,t] \text{ and } x_0 \in B_R \} \quad (12)$$

where $H(\cdot)$ is given in (11b). Obviously, $Y_R(t) \neq \emptyset$ for all $t \geq 0$ and, if (3) holds, the set-valued map $[0,\infty) \ni t \to Y_R(t) \subset \mathbb{R}^{\bar{n}}$ satisfies the CP and further $y(t) \in Y_R(t)$, for every $t \geq t_0 \geq 0$ and $y \in O(t_0, B_R)$; (see notations).

Propositions 2.1 and 2.2 below play a central role to the rest analysis. Parts of their proofs constitute extensions of the methodology employed for the proofs of the main results in Tsinias (2008) and in Boskos and Tsinias (2011). Proposition 2.1 provides sufficient conditions for the solvability of the observer design problem for the case (11) under the additional hypothesis that, it is a priori known that the initial states $x_0$ of (11a) belong to the compact ball $B_R$ of radius $R > 0$ centered at zero:

**Proposition 2.1:** *Consider the system (11) and assume that it is forward complete, namely, (3) holds for certain $\beta \in NN$. We also make the following hypotheses. We assume that there exist a function $g \in C^1(\mathbb{R}_{\geq 0}; \mathbb{R})$ satisfying:*

$$0 < g(t) < 1, \ \forall t \geq 0; \quad (13a)$$
$$\dot{g}(t) \geq -g(t), \ \forall t \geq 0; \quad (13b)$$
$$\lim_{t \to \infty} g(t) = 0 \quad (13c)$$

*an integer $\ell \in \mathbb{N}$, a map $A \in C^0(\mathbb{R}_{\geq 0} \times \mathbb{R}^\ell \times \mathbb{R}^n \times \mathbb{R}^n \times \mathbb{R}^{\bar{n}}; \mathbb{R}^{n \times n})$ and constants $L > 1$, $c_1, c_2 > 0$, $R > 0$ with*

$$c_1 \geq \tfrac{1}{2} \quad (14)$$

*such that the following properties hold:*
**A1.** *For every $\xi \geq 1$ there exists a set-valued map*

$$[0,\infty) \ni t \to Q_R(t) \coloneqq Q_{R,\xi}(t) \subset \mathbb{R}^\ell \quad (15)$$

with $Q_R(t) \neq \emptyset$ for any $t \geq 0$, satisfying the CP and such that
$$\forall t \geq 0, \ y \in Y_R(t), x, z \in \mathbb{R}^n \text{ with } |x| \leq \beta(t, R) \text{ and } |x - z| \leq \xi \Rightarrow$$
$$\Delta F(t, x, z, y) := F(t, x, y) - F(t, z, y)$$
$$= A(t, q, x, x - z, y)(x - z) \text{ for certain } q \in Q_R(t) \tag{16}$$

with $Y_R(\cdot)$ as given by (12);

**A2.** *For every $\xi \geq 1$ there exists a set-valued map $Q_R := Q_{R,\xi}$ as in Hypothesis A1, in such a way that for every $t_0 \geq 0$, a time-varying symmetric matrix $P_R := P_{R,\xi,t_0} \in C^1([t_0, \infty); \mathbb{R}^{n \times n})$ and a function $d_R := d_{R,\xi,t_0} \in C^0([t_0, \infty); \mathbb{R})$ can be found, satisfying:*

$$P_R(t) \geq I_{n \times n}, \ \forall t \geq t_0; \ |P_R(t_0)| \leq L; \tag{17a}$$

$$d_R(t) > c_1, \ \forall t \geq t_0 + 1; \ \int_{t_1}^{t_2} d_R(s) ds > -c_2, \ \forall t_2 \geq t_1 \geq t_0; \tag{17b}$$

$$e' P_R(t) A(t, q, x, e, y) e + \tfrac{1}{2} e' \dot{P}_R(t) e \leq -d_R(t) e' P_R(t) e,$$
$$\forall t \geq t_0, q \in Q_R(t), x \in \mathbb{R}^n, e \in \ker H(t), y \in Y_R(t): \tag{17c}$$
$$|x| \leq \beta(t, R), |e| \leq \xi, e' P_R(t) e \geq g(t)$$

with $Y_R(\cdot)$ as given by (12).
Then, under A1 and A2 the following hold:
(i) For every $\xi \geq 1$ and $t_0 \geq 0$ there exist functions $\bar{d}_R := \bar{d}_{R,\xi,t_0} \in C^0([t_0, \infty); \mathbb{R})$ with

$$\bar{d}_R(t) < d_R(t), \ \forall t \geq t_0; \tag{18a}$$

$$\bar{d}_R(t) \geq c_1, \ \forall t \geq t_0 + 1; \ \int_{t_1}^{t_2} \bar{d}_R(s) ds > -2c_2, \ \forall t_2 \geq t_1 \geq t_0 \tag{18b}$$

and $\phi_R := \phi_{R,\xi,t_0} \in C^1([t_0, \infty); \mathbb{R}_{>0})$ such that

$$e' P_R(t) A(t, q, x, e, y) e + \tfrac{1}{2} e' \dot{P}_R(t) e \leq \phi_R(t) |H(t) e|^2 - \bar{d}_R(t) e' P_R(t) e,$$
$$\forall t \geq t_0, q \in Q_R(t), x \in \mathbb{R}^n, e \in \mathbb{R}^n, y \in Y_R(t): \tag{18c}$$
$$|x| \leq \beta(t, R), |e| \leq \xi, e' P_R(t) e \geq g(t)$$

(ii) Furthermore, the following property holds: For each $\bar{t}_0 \geq t_0 \geq 0$ and constant $\xi$ satisfying

$$\xi \geq \sqrt{L} \exp\{2c_2\} \beta(\bar{t}_0, \bar{R}) \tag{19a}$$
$$\bar{R} := R + 1 \tag{19b}$$

the solution $z(\cdot)$ of system

$$\dot{z} = G_{\bar{t}_0}(t, z, y) := F(t, z, y) + \phi_R(t) P_R^{-1}(t) H'(t)(y - H(t) z) \tag{20a}$$
$$\text{with initial } z(\bar{t}_0) = 0 \tag{20b}$$

where $P_R(\cdot) := P_{R,\xi,\bar{t}_0}(\cdot)$ is given in A2 and with $\phi_R := \phi_{R,\xi,\bar{t}_0} \in C^1([\bar{t}_0, \infty); \mathbb{R}_{>0})$ as given in (18c), is defined for all $t \geq \bar{t}_0$ and the error $e(t) := x(t) - z(t)$ between the trajectory $x(\cdot) := x(\cdot, t_0, x_0)$ of (11a), initiated from $x_0 \in B_R$ at time $t_0 \geq 0$ and the trajectory $z(\cdot) := z(\cdot, \bar{t}_0, 0)$ of (20) satisfies:

- $|e(t)| < \xi$, $\forall t \geq \bar{t}_0$;  (21a)
- $|e(t)| \leq \max\{\xi \exp\{-c_1(t-(\bar{t}_0+1))\}, \sqrt{g(t)}\}$, $\forall t \geq \bar{t}_0 + 1$  (21b)

*It turns out from (13c) and (21b) that for $\bar{t}_0 := t_0$ the ODP is solvable for (11) with respect to $B_R$; particularly the error $e(\cdot)$ between the trajectory $x(\cdot) := x(\cdot, t_0, x_0)$, $x_0 \in B_R$ of (11a) and the trajectory $z(\cdot) := z(\cdot, t_0, z_0)$, $z_0 = 0$ of the observer $\dot{z} = G_{t_0}(t, z, y)$ satisfies (5).*

**Remark:** (i) Notice that (17b) is fulfilled when $d_R(\cdot)$ satisfies $d_R(t) > c_1$, for all $t \geq t_0$.

(ii) The result of Proposition 2.1 is also valid, by considering arbitrary $c_1 > 0$, provided that, instead of (13b), it holds $\dot{g}(t) \geq -2c_1 g(t)$ for all $t \geq 0$.

Based on estimations (21a) and (21b) of Proposition 2.1, we may establish the following proposition which establishes sufficient conditions for the existence of a switching observer exhibiting the state determination of (11), without any a priori information concerning the initial condition.

**Proposition 2.2:** *In addition to the hypothesis of forward completeness for (11a), assume that system (11) satisfies the following property:*

**A3.** *There exist an integer $\ell \in \mathbb{N}$, a map $A \in C^0(\mathbb{R}_{\geq 0} \times \mathbb{R}^\ell \times \mathbb{R}^n \times \mathbb{R}^n \times \mathbb{R}^{\bar{n}}; \mathbb{R}^{n \times n})$, constants $L > 1$, $c_1, c_2 > 0$ such that (14) holds, a function $g \in C^1(\mathbb{R}_{\geq 0}; \mathbb{R})$ satisfying (13) and in such a way that for every $R > 0$ both hypotheses A1 and A2 hold. Then the SODP is solvable for (11) with respect to $\mathbb{R}^n$.*

In order to prove Propositions 2.1 and 2.2 and the main result of Theorem 1.1 in the next section, we first need to establish a preliminary technical result (Lemma 2.1 below). Let $\ell, m, n, \bar{n} \in \mathbb{N}$, consider a pair $(H, A)$ of continuous mappings:

$$H := H(t) : \mathbb{R}_{\geq 0} \to \mathbb{R}^{\bar{n} \times m} ; \quad (22a)$$

$$A := A(t, q, x, e, y) : \mathbb{R}_{\geq 0} \times \mathbb{R}^\ell \times \mathbb{R}^n \times \mathbb{R}^m \times \mathbb{R}^{\bar{n}} \to \mathbb{R}^{m \times m} \quad (22b)$$

and make the following hypothesis:

**H1.** *Let $g(\cdot) \in C^0(\mathbb{R}_{\geq 0}; \mathbb{R})$ satisfying (13a) and assume that for certain constant $R > 0$ and for every $\xi \geq 1$, there exist a function $\beta_R := \beta_{R,\xi} \in N$ and set-valued mappings $[0, \infty) \ni t \to Y_R(t) := Y_{R,\xi}(t) \subset \mathbb{R}^{\bar{n}}$ and $[0, \infty) \ni t \to Q_R(t) := Q_{R,\xi}(t) \subset \mathbb{R}^\ell$ with $Y_R(t) \neq \emptyset$ and $Q_R(t) \neq \emptyset$ for all $t \geq 0$, satisfying the CP, in such a way that for every $t_0 \geq 0$, a time-varying symmetric matrix $P_R := P_{R,\xi,t_0} \in C^1([t_0, \infty); \mathbb{R}^{m \times m})$ can be found satisfying*

$$P_R(t) \geq I_{m \times m}, \quad \forall t \geq t_0 \quad (23)$$

*and a function $d_R := d_{R,\xi,t_0} \in C^0([t_0, \infty); \mathbb{R})$, in such a way that*

$$\begin{aligned}&e'P_R(t)A(t,q,x,e,y)e+\tfrac{1}{2}e'\dot{P}_R(t)e\leq -d_R(t)e'P_R(t)e,\\&\forall t\geq t_0,\,q\in Q_R(t),\,x\in\mathbb{R}^n,\,e\in\ker H(t),\,y\in Y_R(t):\\&|x|\leq \beta_R(t),\,|e|\leq \xi,\,e'P_R(t)e\geq g(t)\end{aligned} \qquad (24)$$

**Lemma 2.1:** *Consider the pair $(H,A)$ of the continuous time-varying mappings in (22) and assume that H1 is fulfilled for certain $R>0$. Then for every $\xi\geq 1$, $t_0\geq 0$ and $\bar{d}_R:=\bar{d}_{R,\xi,t_0}\in C^0([t_0,\infty);\mathbb{R})$ with*

$$\bar{d}_R(t)<d_R(t),\;\forall t\geq t_0 \qquad (25)$$

*there exists a function $\phi_R:=\phi_{R,\xi,t_0}\in C^1([t_0,\infty);\mathbb{R}_{>0})$ such that*

$$\begin{aligned}&e'P_R(t)A(t,q,x,e,y)e+\tfrac{1}{2}e'\dot{P}_R(t)e\leq \phi_R(t)|H(t)e|^2-\bar{d}_R(t)e'P_R(t)e,\\&\forall t\geq t_0,\,q\in Q_R(t),\,x\in\mathbb{R}^n,\,e\in\mathbb{R}^m,\,y\in Y_R(t):\\&|x|\leq \beta_R(t),\,|e|\leq \xi,\,e'P_R(t)e\geq g(t)\end{aligned} \qquad (26)$$

**Proof:** Let $R$ as given in our statement, let $\xi\geq 1$ and $t_0\geq 0$ and consider the mappings $Y_R$, $Q_R$, $P_R$ and $d_R$ as involved in H1 and a function $\bar{d}_R$ satisfying (25). Also, consider the following continuous function $\rho_R:=\rho_{R,\xi,t_0}:[t_0,\infty)\times\{w\in\mathbb{R}^m:|w|=1\}\to\mathbb{R}$:

$$\rho_R(t,w):=\sqrt{\frac{g(t)}{w'P_R(t)w}} \qquad (27)$$

By taking into account (13a),(23), definition (27) and recalling that $\xi\geq 1$, it follows that $\rho_R(t,w)<\xi$, for every $t\in[t_0,\infty)$ and $w\in\mathbb{R}^m$ with $|w|=1$. Then for each $(t,w)\in[t_0,\infty)\times\{w\in\mathbb{R}^m:|w|=1\}$ we may consider the following nonempty closed interval $I_R:=I_{R,\xi,t_0}\subset\mathbb{R}$ defined as:

$$I_R(t,w):=[\rho_R(t,w),\xi] \qquad (28)$$

Obviously, the set-valued map $[t_0,\infty)\times\{w\in\mathbb{R}^m:|w|=1\}\ni(t,w)\to I_R(t,w)$ satisfies the CP. Also, due to (27) and (28), it follows that for any subspace $X$ of $\mathbb{R}^m$ we have:

$$\{e\in X:|e|\leq \xi\text{ and }e'P_R(t)e\geq g(t)\}=\{w\sigma:w\in X\text{ and }|w|=1,\sigma\in I_R(t,w)\} \qquad (29)$$

Using (29) with $X:=\ker H(t)$, it follows that (24) is equivalent to:

$$\begin{aligned}&w'P_R(t)A(t,q,x,w\sigma,y)w+\tfrac{1}{2}w'\dot{P}_R(t)w\leq -d_R(t)w'P_R(t)w,\\&\forall t\geq t_0,\,q\in Q_R(t),\,x\in\mathbb{R}^n,\,w\in\ker H(t),\,\sigma\in I_R(t,w),\,y\in Y_R(t):\\&|x|\leq \beta_R(t),\,|w|=1\end{aligned} \qquad (30)$$

Also, by applying (29) with $X:=\mathbb{R}^m$, the desired (26) is equivalent to the existence of a function $\phi_R\in C^1([t_0,\infty);\mathbb{R}_{>0})$ such that the following implication holds:

$$w \in \mathbb{R}^m : |w|=1 \Rightarrow w'P_R(t)A(t,q,x,w\sigma,y)w + \tfrac{1}{2}w'\dot{P}_R(t)w \leq \phi_R(t)|H(t)w|^2 - \bar{d}_R(t)w'P_R(t)w,$$
$$\forall t \geq t_0, \, q \in Q_R(t), \, |x| \leq \beta_R(t), \, \sigma \in I_R(t,w), \, y \in Y_R(t)$$
(31)

In order to show (31), we proceed as follows. Consider the map $D_R := D_{R,\xi,t_0} : [t_0,\infty) \times \mathbb{R}^\ell \times \mathbb{R}^n \times \mathbb{R}^m \times \mathbb{R} \times \mathbb{R}^{\bar{n}} \to \mathbb{R}$ and the set $K_R(t) := K_{R,\xi,t_0}(t)$, $t \geq t_0$, defined as:

$$D_R(t,q,x,e,\sigma,y) := e'P_R(t)A(t,q,x,e\sigma,y)e + \tfrac{1}{2}e'\dot{P}_R(t)e + \bar{d}_R(t)e'P_R(t)e; \quad (32a)$$

$$K_R(t) := \{|w|=1 : D_R(t,q,x,w,\sigma,y) < 0, \text{ for every } q \in Q_R(t), \, |x| \leq \beta_R(t), \, \sigma \in I_R(t,w), \, y \in Y_R(t)\} \quad (32b)$$

Notice that for those $t \geq t_0$ for which $\operatorname{rank} H(t) < m$, the set $K_R(t)$ is nonempty, since we have:

$$w \in \ker H(t) \text{ and } |w|=1 \Rightarrow w \in K_R(t) \quad (33)$$

Indeed, let $w \in \ker H(t)$ for certain $w \in \mathbb{R}^m$ with $|w|=1$. Then, by using (23),(25),(30) and (32a), we deduce $D_R(t,q,x,w,\sigma,y) \leq (\bar{d}_R(t) - d_R(t))w'P_R(t)w < 0$, for all $q \in Q_R(t)$, $|x| \leq \beta_R(t)$, $\sigma \in I_R(t,w)$ and $y \in Y_R(t)$, therefore, (32b) asserts that $w \in K_R(t)$ and this establishes (33). Consequently, $K_R(t) \neq \emptyset$, provided that $\operatorname{rank} H(t) < m$. For each $t \geq t_0$, define $K_R^c(t) := \{w \in \mathbb{R}^m : |w|=1, w \notin K_R(t)\}$. Then by (32b) it follows:

$$K_R^c(t) = \{|w|=1 : D_R(t,q,x,w,\sigma,y) \geq 0,$$
$$\text{for some } q \in Q_R(t), |x| \leq \beta_R(t), \sigma \in I_R(t,w), y \in Y_R(t)\} \quad (34)$$

We now prove that for every $t \geq t_0$ the set $K_R^c(t)$ is closed. Indeed, let $t \geq t_0$ and without any loss of generality assume that $K_R^c(t) \neq \emptyset$. We show that for every sequence $(w_\nu)_{\nu \in \mathbb{N}} \subset K_R^c(t)$ with $w_\nu \to w$ we have $w \in K_R^c(t)$. Define:

$$\bar{Q}_R(t,w) := Q_R(t) \times \{x \in \mathbb{R}^n : |x| \leq \beta_R(t)\} \times I_R(t,w) \times Y_R(t), \, (t,w) \in [t_0,\infty) \times \{w \in \mathbb{R}^m : |w|=1\}$$
(35)

Then, since $w_\nu \in K_R^c(t)$, it follows from (34) and (35) that there exists $(q_\nu, x_\nu, \sigma_\nu, y_\nu) \in \bar{Q}_R(t,w_\nu)$ with $D_R(t,q_\nu,x_\nu,w_\nu,\sigma_\nu,y_\nu) \geq 0$. Notice that the set valued map $[t_0,\infty) \times \{w \in \mathbb{R}^m : |w|=1\} \ni (t,w) \to \bar{Q}_R(t,w) \subset \mathbb{R}^\ell \times \mathbb{R}^n \times \mathbb{R} \times \mathbb{R}^{\bar{n}}$ satisfies the CP, hence, continuity of the map $D_R$ implies that $\underline{\lim} D_R(t,q_\nu,x_\nu,w_\nu,\sigma_\nu,y_\nu) \to D_R(t,q,x,w,\sigma,y) \geq 0$, for certain $(q,x,\sigma,y) \in \bar{Q}_R(t,w)$. It follows from (34) and (35) that $w \in K_R^c(t)$, therefore $K_R^c(t)$ is closed. Next, consider the following map $\omega_R := \omega_{R,\xi,t_0} : [t_0,\infty) \to [0,\infty]$:

$$\omega_R(t) := \begin{cases} \min\{|H(t)w| : w \in K_R^c(t)\}, & \text{if } K_R^c(t) \neq \emptyset \\ \infty & \text{if } K_R^c(t) = \emptyset \end{cases} \quad (36)$$

Notice that, if $K_R^c(t) \neq \emptyset$, then $\{|H(t)w| : w \in K_R^c(t)\}$ is nonempty and compact, thus $\omega_R(\cdot)$ is well defined. Consequently, by (33) and (36), we have $\omega_R(t) > 0$ for every $t \geq t_0$. Moreover, it holds:

$$\inf\{\omega_R(t): t \in [t_0, T]\} > 0, \; \forall T > t_0 \tag{37}$$

Indeed, suppose on the contrary that $\omega_R(t_\nu) \to 0$ for certain sequence $(t_\nu)_{\nu \in \mathbb{N}} \subset [t_0, T]$ with $t_\nu \to t$. Then by (36) we may assume that without any loss of generality $K_R^c(t_\nu) \neq \emptyset$ for every $\nu \in \mathbb{N}$, hence, there exists $w_\nu \in K_R^c(t_\nu)$ such that $|H(t_\nu)w_\nu| \to 0$. Since $|w_\nu| = 1$, we may also assume that, without any loss of generality, $w_\nu \to w$ for some $w \in \mathbb{R}^m$ with $|w| = 1$. Therefore, continuity of $H(\cdot)$ implies $H(t)w = 0$, hence, by (33) we have that $w \in K_R(t)$. On the other hand $w_\nu \in K_R^c(t_\nu)$, which by virtue of (34) and (35), implies that $D_R(t_\nu, q_\nu, x_\nu, w_\nu, \sigma_\nu, y_\nu) \geq 0$ for some $(q_\nu, x_\nu, \sigma_\nu, y_\nu) \in \bar{Q}_R(t_\nu, w_\nu)$. The latter in conjunction with the compactness property for $\bar{Q}_R(\cdot, \cdot)$ and continuity of the map $D_R$ implies that $\underline{\lim} D_R(t_\nu, q_\nu, x_\nu, w_\nu, \sigma_\nu, y_\nu) \to D_R(t, q, x, w, \sigma, y) \geq 0$ for some $(q, x, \sigma, y) \in \bar{Q}_R(t, w)$, therefore by (34) and (35) we get $w \in K_R^c(t)$, which is a contradiction. This establishes (37). Consider now the function $\bar{\omega}_R := \bar{\omega}_{R, \xi, t_0} : [t_0, \infty) \to [0, \infty)$ defined as:

$$\bar{\omega}_R(t) := \begin{cases} \dfrac{1}{\omega_R^2(t)}, & K_R^c(t) \neq \emptyset \\ 0, & K_R^c(t) = \emptyset \end{cases} \tag{38}$$

By taking into account (37),(38) and recalling again the CP for $Q_R(\cdot)$ and $Y_R(\cdot)$, it follows that there exists a function $\phi_R := \phi_{R, \xi, t_0} \in C^1([t_0, \infty); \mathbb{R}_{>0})$ which satisfies:

$$\begin{aligned} \phi_R(t) > \bar{\omega}_R(t) \sup\{|P_R(t)||A(t, q, x, e, y)| + \tfrac{1}{2}|\dot{P}_R(t)| + |\bar{d}_R(t)||P_R(t)|: \\ q \in Q_R(t), |x| \leq \beta_R(t), |e| \leq \xi, y \in Y_R(t)\}, \forall t \geq t_0 \end{aligned} \tag{39}$$

We are now in a position to show that (31) holds with $\phi_R(\cdot)$ as above. We distinguish two cases. First, consider those $t \geq t_0$ for which $K_R^c(t) \neq \emptyset$. Then, if $K_R(t) \neq \emptyset$ and $w \in K_R(t)$, implication (31) is a direct consequence of (32) and the fact that $\phi_R(t) > 0$. When $w \in K_R^c(t)$, in order to show implication (31), it suffices by virtue of (27),(28) and (36) to prove:

$$\begin{aligned} \sup\{|P_R(t)||A(t, q, x, w\sigma, y)| + \tfrac{1}{2}|\dot{P}_R(t)| + |\bar{d}_R(t)||P_R(t)|: \\ q \in Q_R(t), |x| \leq \beta_R(t), |w| = 1, 0 \leq \sigma \leq \xi, y \in Y_R(t)\} \leq \phi_R(t)\omega_R^2(t) \end{aligned}$$

The latter is a consequence of (38) and (39). Now, consider those $t \geq t_0$ for which $K_R^c(t) = \emptyset$, or equvalently $K_R(t) = \{w \in \mathbb{R}^m : |w| = 1\}$. Then (31) is a consequence of (32). ∎

We are now in a position to prove Proposition 2.1.

**Proof of Proposition 2.1:** The proof of Statement (i) is based on the result of Lemma 2.1. Let $R > 0$ as given in our statement and let $\xi \geq 1$ and $t_0 \geq 0$. Consider a function $d_R := d_{R, \xi, t_0}$ as determined in Hypothesis A2 and let $c := \min\{d_R(t_0 + 1) - c_1, c_2\}$, which due to the first inequality of (17b), is strictly positive. Then, by taking into account the first inequality of (17b), we can determine a function $\bar{d}_R := \bar{d}_{R, \xi, t_0} \in C^0([t_0, \infty); \mathbb{R})$ satisfying:

$$\bar{d}_R(t) = d_R(t) - c, \ \forall t \in [t_0, t_0 + 1]; \tag{40a}$$

$$c_1 \leq \bar{d}_R(t) < d_R(t), \ \forall t \geq t_0 + 1 \tag{40b}$$

The desired (18a) and the first inequality of (18b) are both direct consequences of (40). Also, from the definition of the constant $c$ above and by taking into account the second inequality of (17b) and (40a), we have $\int_{t_1}^{t_2} \bar{d}_R(s) ds = \int_{t_1}^{t_2} (d_R(s) - c) ds = \int_{t_1}^{t_2} d_R(s) ds - c(t_2 - t_1) > -2c_2$, for every $t_2 \geq t_1$, $t_1, t_2 \in [t_0, t_0 + 1]$, which, due to (40b) and the fact that $c_1 > 0$, asserts that the second inequality in (18b) is fulfilled. In order to show (18c), consider the pair $(H, A)$ as given in (22) with $m = n$ and where $\ell, n, \bar{n}$, $H(t)$ and $A(t, q, x, e, y)$ are involved in Hypotheses A1 and A2. Notice that, due to A2 and particularly by taking into account (17c), first inequality of (17a) and the fact that $Y_R(\cdot)$ satisfies the CP, it follows that Property H1 holds for the pair $(H, A)$ with $Y_R(\cdot)$, $Q_R(\cdot)$, $P_R(\cdot)$, $d_R(\cdot)$, $g(\cdot)$ and $\beta_R(\cdot) := \beta(\cdot, R)$ as given in (12),(15),(17a),(17b),(13) and (3), respectively. Therefore, by taking into account that the function $\bar{d}_R$ above satisfies (18a) and by invoking the result of Lemma 2.1, it follows that, for the given $R > 0$, $\xi \geq 1$ and $t_0 \geq 0$, a function $\phi_R := \phi_{R,\xi,t_0} \in C^1([t_0, \infty); \mathbb{R}_{>0})$ can be found satisfying (18c) and the proof of Statement (i) is completed. We next proceed to the establishment of Statement (ii), namely, we show that (21a) and (21b) hold. Consider the constant $R > 0$ as given in our statement, let $\bar{t}_0 \geq t_0 \geq 0$, and let $\xi$ be a constant satisfying (19). Notice that assumption (3) guarantees that

$$\beta(t, s) \geq s \text{ for all } t \geq 0 \text{ and } s \geq 0 \tag{41}$$

and, since $L > 1$, $c_2 > 0$ it follows from (19) and (41) that $\xi > 1$. Therefore, we may apply the result of Statement (i) with $\xi$ as above and $t_0 := \bar{t}_0$, in order to find functions $\bar{d}_R \in C^0([\bar{t}_0, \infty); \mathbb{R})$ and $\phi_R \in C^1([\bar{t}_0, \infty); \mathbb{R}_{>0})$ satisfying (18). Let $e(\cdot) := x(\cdot, t_0, x_0) - z(\cdot, \bar{t}_0, z_0)$ be the error between the trajectory of (11a) and the trajectory of (20a) with $x_0 \in B_R$ and $z_0 = z(\bar{t}_0) = 0$. Then, by invoking (11b), the error $e(\cdot)$ satisfies:

$$\begin{aligned}\dot{e} &= F(t, x, y) - F(t, z, y) - \phi_R(t) P_R^{-1}(t) H'(t)(y - H(t)z) \\ &= F(t, x, y) - F(t, z, y) - \phi_R(t) P_R^{-1}(t) H'(t) H(t) e, \ t \geq \bar{t}_0\end{aligned} \tag{42a}$$

and obviously (due to (20b)):

$$e(\bar{t}_0) = x(\bar{t}_0, t_0, x_0) - z(\bar{t}_0, \bar{t}_0, z_0) = x(\bar{t}_0, t_0, x_0) \tag{42b}$$

Notice that, since $L > 1$, $c_2 > 0$ and due to (3),(19) and (42b), we have:

$$|e(\bar{t}_0)| \leq \beta(\bar{t}_0, R) < \xi, \ \forall x_0 \in B_R \tag{43}$$

**Establishment of (21a):** We show that for each $x_0 \in B_R$ it holds that $|e(t)| < \xi$ for all $t \in [\bar{t}_0, T_{\max})$, and therefore $T_{\max} = \infty$, where $T_{\max}$ denotes the maximum existence time of the solution $e(\cdot)$ of (42a) corresponding to certain $x_0 \in B_R$ and satisfying (42b). Suppose on the contrary, by taking into account (43), that there would exist $\bar{t} \in (\bar{t}_0, T_{max})$ and $x_0 \in B_R$ such that

$$|e(\bar{t})| = \xi ; \qquad (44a)$$

$$|e(t)| < \xi, \ \forall t \in [\bar{t}_0, \bar{t}) \qquad (44b)$$

By invoking (12) it follows that, since $x_0 \in B_R$, the output $y(\cdot)$ of (11) satisfies:

$$y(t) \in Y_R(t), \ \forall t \in [\bar{t}_0, \infty) \qquad (45)$$

and the latter in conjunction with (16), together with (3) and (44), implies:

$$\Delta F(t, x(t), z(t), y(t)) = A(t, q, x(t), x(t) - z(t), y(t))(x(t) - z(t)), \forall t \in [\bar{t}_0, \bar{t}] \qquad (46)$$
$$\text{for some } q := q(t) \in Q_R(t)$$

Define:

$$V(t, e) := e' P_R(t) e, \ t \ge \bar{t}_0, \ e \in \mathbb{R}^n \qquad (47)$$

and evaluate its time-derivative along the trajectories $e(\cdot)$ of (42a). By (46) and (47) it follows:

$$\dot{V}(t, e(t)) = e'(t) \dot{P}_R(t) e(t) + 2 e'(t) P_R(t) A(t, q, x(t), e(t), y(t)) e(t) - 2 \phi_R(t) |H(t) e(t)|^2, \qquad (48)$$
$$\forall t \in [\bar{t}_0, \bar{t}] \text{ and for some } q := q(t) \in Q_R(t)$$

Also, notice that, due to the second inequality of (17a),(43) and definition (47), we have:

$$V(\bar{t}_0, e(\bar{t}_0)) = e'(\bar{t}_0) P_R(\bar{t}_0) e(\bar{t}_0) \le L |e(\bar{t}_0)|^2 \le L \beta^2(\bar{t}_0, \bar{R}) \qquad (49)$$

**Claim 1:** We claim that

$$V(t, e(t)) < L \exp\{4 c_2\} \beta^2(\bar{t}_0, \bar{R}), \ \forall t \in [\bar{t}_0, \bar{t}] \qquad (50)$$

**Proof of Claim 1:** Suppose on the contrary, by taking into account (49), that there would exist $\tau \in (\bar{t}_0, \bar{t}]$ such that

$$V(\tau, e(\tau)) \ge L \exp\{4 c_2\} \beta^2(\bar{t}_0, \bar{R}) > L \beta^2(\bar{t}_0, \bar{R}) \qquad (51)$$

where, in order to obtain the second inequality, we have taken into account that $L > 0$ and $c_2 > 0$. Then by (13a),(49) and (51) there would exist $\bar{\tau} \in [\bar{t}_0, \tau)$, such that

$$V(\bar{\tau}, e(\bar{\tau})) = L \beta^2(\bar{t}_0, \bar{R}) ; \ V(t, e(t)) \ge L \beta^2(\bar{t}_0, \bar{R}) \ge 1 \ge g(t), \ \forall t \in [\bar{\tau}, \tau] \qquad (52)$$

where, the second inequality is a consequence of the choices $L > 1$, $\bar{R} > 1$ (see (19)) and inequality (41). By invoking our assumption (18c), together with (3),(44),(45),(47),(48),(52) and recalling that $x_0 \in B_R$, it follows:

$$\dot{V}(t, e(t)) \le -2 \bar{d}_R(t) V(t, e(t)), \ \forall t \in [\bar{\tau}, \tau] \qquad (53)$$

From (53) we get $V(t,e(t)) \leq V(\bar{\tau},e(\bar{\tau}))\exp\left(-2\int_{\bar{\tau}}^{t}\bar{d}_R(s)ds\right)$ for all $t \in [\bar{\tau},\tau]$, therefore, by recalling the second inequality of (18b) and (52), it follows that $V(\tau,e(\tau)) \leq V(\bar{\tau},e(\bar{\tau}))\exp\left(-2\int_{\bar{\tau}}^{\tau}\bar{d}_R(s)ds\right) < L\exp\{4c_2\}\beta^2(\bar{t}_0,\bar{R})$. The latter contradicts (51), thus (50) is satisfied and the proof of Claim 1 is completed. ◄

By invoking the first inequality of (17a), it follows from (19a),(47) and (50) that $|e(t)| < \sqrt{L}\exp\{2c_2\}\beta(\bar{t}_0,\bar{R}) \leq \xi$ for all $t \in [\bar{t}_0,\bar{t}]$, thus $|e(\bar{t})| < \xi$, which in conjunction with (44a) implies that $T_{\max} = \infty$, hence, (21a) is fulfilled. The latter in conjunction with forward completeness of (11a) asserts that the solution $z(\cdot) := z(\cdot,\bar{t}_0,0)$ of (20) is defined for all $t \geq \bar{t}_0$.

**Establishment of (21b):** First, we notice that, since $T_{\max} = \infty$, and, due to (21a), relations (48) and (50) hold for all $t \geq \bar{t}_0$ and $x_0 \in B_R$, namely:

$$\dot{V}(t,e(t)) = e'(t)\dot{P}_R(t)e(t) + 2e'(t)P_R(t)A(t,q,x(t),e(t),y(t))e(t) - 2\phi_R(t)|H(t)e(t)|^2,$$
$$\forall t \geq \bar{t}_0 \text{ and for some } q := q(t) \in Q_R(t)$$
(54a)

$$V(t,e(t)) < L\exp\{4c_2\}\beta^2(\bar{t}_0,\bar{R}), \quad \forall t \geq \bar{t}_0$$
(54b)

To complete the proof of (21b), we need to establish the following claim:

**Claim 2:** For every $x_0 \in B_R$ the following implication holds:

$$V(\bar{t},e(\bar{t})) \leq g(\bar{t}) \text{ for certain } \bar{t} \geq \bar{t}_0 + 1 \Rightarrow V(t,e(t)) \leq g(t) \text{ for all } t > \bar{t} \quad (55)$$

**Proof of Claim 2:** Indeed, suppose on the contrary that there would exist $x_0 \in B_R$ such that

$$V(\bar{t},e(\bar{t})) \leq g(\bar{t}) \text{ for certain } \bar{t} \geq \bar{t}_0 + 1 \text{ and } V(\tau,e(\tau)) > g(\tau) \text{ for some } \tau > \bar{t} \quad (56)$$

Then by (56), there would exist $\bar{\tau} \in [\bar{t},\tau)$ such that

$$V(\bar{\tau},e(\bar{\tau})) = g(\bar{\tau}); \ V(t,e(t)) > g(t), \ \forall t \in (\bar{\tau},\tau] \quad (57)$$

By invoking assumption (18c) together with (3),(21a),(45),(47),(54a),(57) and recalling that $x_0 \in B_R$, we deduce $\dot{V}(t,e(t)) \leq -2\bar{d}_R(t)V(t,e(t))$ for all $t \in [\bar{\tau},\tau]$, thus, by recalling that $\bar{\tau} \geq \bar{t}_0 + 1$ and applying the first inequality of (18b), we have:

$$\dot{V}(t,e(t)) \leq -2c_1 V(t,e(t)), \ \forall t \in [\bar{\tau},\tau] \quad (58)$$

Define $h(t) := V(t,e(t)) - g(t)$, $t \in [\bar{\tau},\tau]$. Then (57) is equivalently written: $h(\bar{\tau}) = 0$ and $h(t) > 0$ for all $t \in (\bar{\tau},\tau]$, thus there would exist $T \in (\bar{\tau},\tau)$ such that the time derivative $\dot{h}$ of $h$ satisfies $\dot{h}(T) = \dot{V}(T,e(T)) - \dot{g}(T) > 0$. The latter in conjunction with (13b) and (58) implies $c_1 V(T,e(T)) < \frac{1}{2}g(T)$, hence, by (14) we get $V(T,e(T)) < g(T)$. On the other hand, it follows from (57) that $V(T,e(T)) > g(T)$, a contradiction, therefore (55) is established. ◄

We are now in position to complete the proof of (21b). Let $x_0 \in B_R$ and consider the following cases:

**Case 1:**
$$V(\overline{t}_0+1, e(\overline{t}_0+1)) \le g(\overline{t}_0+1) \tag{59}$$

Then, by taking into account (55) and (59), it follows that $V(t,e(t)) \le g(t)$ for all $t \ge \overline{t}_0+1$, therefore, by invoking the first inequality of (17a) and (47), we get $|e(t)| \le \sqrt{g(t)}$ for every $t \ge \overline{t}_0+1$, which implies (21b).

**Case 2:**
$$V(\overline{t}_0+1, e(\overline{t}_0+1)) > g(\overline{t}_0+1) \text{ and } V(t,e(t)) > g(t), \forall t > \overline{t}_0+1 \tag{60}$$

Then, by invoking (3),(18c),(21a),(45),(47),(54a),(60) and the fact that $x_0 \in B_R$, we deduce that $\dot{V}(t,e(t)) \le -2\overline{d}_R(t)V(t,e(t))$ for all $t \ge \overline{t}_0+1$. It then follows, by invoking the first inequality of (18b), that $\dot{V}(t,e(t)) \le -2c_1 V(t,e(t))$ for every $t \ge \overline{t}_0+1$, therefore, by (54b) we get:

$$V(t,e(t)) \le V(\overline{t}_0+1, e(\overline{t}_0+1))\exp\{-2c_1(t-(\overline{t}_0+1))\}$$
$$< L\exp\{4c_2\}\beta^2(\overline{t}_0, \overline{R})\exp\{-2c_1(t-(\overline{t}_0+1))\}, \forall t \ge \overline{t}_0+1$$

The latter in conjunction with (17a),(19a) and (47) yields $|e(t)| < \xi\exp\{-c_1(t-(\overline{t}_0+1))\}$ for all $t \ge \overline{t}_0+1$, hence, (21b) is fulfilled.

**Case 3:**
$$V(\overline{t}_0+1, e(\overline{t}_0+1)) > g(\overline{t}_0+1) \text{ and there exists } \overline{t} > \overline{t}_0+1 \text{ such that } V(\overline{t}, e(\overline{t})) \le g(\overline{t}) \tag{61}$$

If (61) holds, then we may define:
$$\tau := \min\{t > \overline{t}_0+1 : V(t,e(t)) = g(t))\} \tag{62}$$

It follows from (61) and (62) that $V(t,e(t)) \ge g(t)$, $\forall t \in [\overline{t}_0+1, \tau]$, therefore by arguing as in Case 2 above, we deduce:

$$|e(t)| < \xi\exp\{-c_1(t-(\overline{t}_0+1))\}, \forall t \in [\overline{t}_0+1, \tau] \tag{63}$$

Also, from (62) we have $V(\tau, e(\tau)) = g(\tau)$, therefore, from (17a),(47) and (55) it holds:
$$V(t,e(t)) \le g(t), \forall t \in [\tau, \infty) \Rightarrow |e(t)| \le \sqrt{g(t)}, \forall t \in [\tau, \infty) \tag{64}$$

It follows from (63) and (64) that (21b) is again fulfilled and this completes the proof of Proposition 2.1. ∎

**Proof of Proposition 2.2:** Consider the system (11) initiated at $t_0 \ge 0$. We proceed to the construction of a sequence of times $(t_k)_{k \in \mathbb{N}}$ and an appropriate sequence of mappings $(G_k)_{k \in \mathbb{N}}$, in such a way that the sequence of systems (8) exhibits the state determination of (11), according to the definition of solvability of the SODP given in Section 1. Let $y \in O(t_0, \mathbb{R}^n)$ be an output of (11). Also, let $L > 1$ and $c_1, c_2 > 0$ be constants, such that (14) holds and in such a way that for every $R > 0$ Properties A1 and A2 hold, according to Assumption A3.

**Claim 1 (Induction Hypothesis):** Let $Y_R(\cdot)$ as given by (12), with $R := k = 1, 2, \ldots$, and let $\beta(\cdot, \cdot)$, $g(\cdot)$ as given by (3),(13) respectively. Then, for $L$, $c_1$, $c_2$ as above and for every $k \in \mathbb{N}$ there exist positive constants $\xi_k$ and $t_k$, with $t_1 := t_0$ ($t_0$ being the initial time), a set-valued mapping $Q_k$ satisfying the CP, a time-varying symmetric matrix $P_k \in C^1([t_{k-1}, \infty); \mathbb{R}^{n \times n})$ and functions $\bar{d}_k \in C^0([t_{k-1}, \infty); \mathbb{R})$ and $\phi_k \in C^1([t_{k-1}, \infty); \mathbb{R}_{>0})$ such that

$$\xi_k \geq 1; \tag{65}$$

$$t_1 := t_0;\ t_{k+1} \geq t_k + 1,\ k = 1, 2, \ldots \tag{66a}$$

$$\text{and therefore } \lim_{k \to \infty} t_k = \infty; \tag{66b}$$

$$P_k(t) \geq I_{n \times n},\ \forall t \geq t_{k-1};\ |P_k(t_{k-1})| \leq L; \tag{67}$$

$$\bar{d}_k(t) \geq c_1,\ \forall t \geq t_{k-1} + 1;\ \int_{\tau_1}^{\tau_2} \bar{d}_k(s)ds > -2c_2,\ \forall \tau_2 \geq \tau_1 \geq t_{k-1}; \tag{68}$$

$$e'P_k(t)A(t, q, x, e, y)e + \tfrac{1}{2}e'\dot{P}_k(t)e - \phi_k(t)|H(t)e|^2 \leq -\bar{d}_k(t)e'P_k(t)e,$$
$$\forall t \geq t_{k-1},\ q \in Q_k(t),\ x \in \mathbb{R}^n,\ e \in \mathbb{R}^n,\ y \in Y_k(t): \tag{69}$$
$$|x| \leq \beta(t, k),\ |e| \leq \xi_k,\ e'P_k(t)e \geq g(t)$$

Specifically, for each $k \in \mathbb{N}$ the desired constants $\xi_k$ and $t_k$ above are defined recursively as follows:

- $$\xi_i := \sqrt{L} \exp\{2c_2\}\beta(t_{i-1}, i+1),\ i = 1, \ldots, k \tag{70}$$

- $$t_i := \min\{\tau \geq t_{i-1} + 1 : \max\{\xi_i \exp\{-c_1(t - (t_{i-1} + 1))\}, \sqrt{g(t)}\} \leq \tfrac{1}{i},\ \text{for all } t \geq \tau\}$$
$$i = 2, \ldots, k \text{ for } k \geq 2 \text{ and } t_1 := t_0 \tag{71}$$

**Proof of Claim 1 for $k := 1$:** Define $\xi_1$ as in (70) with $i = 1$, namely:

$$\xi_1 := \sqrt{L} \exp\{2c_2\}\beta(t_0, 2) \tag{72}$$

Since $L > 1$ and $c_2 > 0$ it follows from (41) and (72) that $\xi_1 \geq 1$, namely, (65) holds for $k = 1$. Then according to our assumption that A2 holds for $R = 1$, it follows that for the given $\xi_1$ above, there exist a set-valued map $Q_1 := Q_{1, \xi_1}$ satisfying the CP and a pair of mappings $d_1 \in C^0([t_0, \infty); \mathbb{R})$ and $P_1 \in C^1([t_0, \infty); \mathbb{R}^{n \times n})$, where $P_1(\cdot)$ satisfies (67) and in such a way that for $R = 1$ conditions (17b,c) are fulfilled. Then, by applying the result of Proposition 2.1(i) with $R = 1$ for the given $\xi_1$ and $t_0$ it follows that there exist functions $\bar{d}_1 \in C^0([t_0, \infty); \mathbb{R})$ and $\phi_1 \in C^1([t_0, \infty); \mathbb{R}_{>0})$ satisfying (68) and (69) with $k = 1$. Finally, define $t_1 := t_0$ and the above discussion establishes Claim 1 for the case $k = 1$.

**Proof of Claim 1 (general step of induction procedure):** Suppose next that all requirements of Claim 1 are valid for some $k \in \mathbb{N}$. We show that it also holds for $k := k + 1$. Define $\xi_{k+1}$ as in (70) with $i = k + 1$, namely:

$$\xi_{k+1} := \sqrt{L} \exp\{2c_2\}\beta(t_k, k+2) \tag{73}$$

Since $L > 1$ and $c_2 > 0$ it follows from (41) and (73) that $\xi_{k+1} \geq 1$, thus (65) is fulfilled for $k = k + 1$. Also, for $k = k + 1$ we define the corresponding $t_{k+1}$, $Q_{k+1}(\cdot)$, $P_{k+1}(\cdot)$, $\bar{d}_{k+1}(\cdot)$ and $\phi_{k+1}(\cdot)$ as follows. By again recalling Assumption A3, we may apply A2 with $R = k+1$. It follows that, for the constant $\xi_{k+1}$ as defined by (73), there exist a set-valued map $Q_{k+1} := Q_{k+1,\xi_{k+1}}$ satisfying the CP and a pair of mappings $d_{k+1} \in C^0([t_k,\infty);\mathbb{R})$ and $P_{k+1} \in C^1([t_k,\infty);\mathbb{R}^{n\times n})$ satisfying (67) as well as (17b,c) for $R = k+1$. Then, by again invoking the result of Proposition 2.1(i) with $R = k+1$, for the constant $\xi_{k+1}$ as given by (73), there exist functions $\bar{d}_{k+1} \in C^0([t_k,\infty);\mathbb{R})$ and $\phi_{k+1} \in C^1([t_k,\infty);\mathbb{R}_{>0})$, which satisfy (68) and (69) with $k = k+1$. Finally, define $t_{k+1}$ as in (71) with $i = k+1$, namely:

$$t_{k+1} := \min\{\tau \geq t_k + 1 : \max\{\xi_{k+1}\exp\{-c_1(t-(t_k+1))\},\sqrt{g(t)}\} \leq \tfrac{1}{k+1}, \text{ for all } t \geq \tau\}$$

Since $c_1 > 0$, it can be easily verified by exploiting (13c), that $t_{k+1}$ is well defined and satisfies (66a). This completes the establishment of the general step of the induction hypothesis. ◄

**Construction of the switching observer:** Using the properties of Claim 1, we proceed to an explicit construction of a switching observer for (11). We define the sequence of systems:

$$\dot{z}_k = G_k(t,z_k,y), \quad t \in [t_{k-1},t_{k+1}], \tag{74a}$$

$$\text{with initial } z_k(t_{k-1}) = 0, \quad k = 1,2,\ldots \tag{74b}$$

$$G_k(t,z,y) \begin{cases} := F(t,z,y) + \phi_k(t)P_k^{-1}(t)H'(t)(y - H(t)z) \\ \quad \text{for } |z| \leq \zeta_k, t \in [t_{k-1},t_{k+1}], \\[4pt] := (F(t,z,y) + \phi_k(t)P_k^{-1}(t)H'(t)(y - H(t)z))\dfrac{2\zeta_k - |z|}{\zeta_k} \\ \quad \text{for } \zeta_k \leq |z| \leq 2\zeta_k, t \in [t_{k-1},t_{k+1}], \\[4pt] := 0, \quad \text{for } |z| \geq 2\zeta_k, t \in [t_{k-1},t_{k+1}] \end{cases} \tag{74c}$$

$$\zeta_k := \beta(t_{k+1},k) + \xi_k \tag{74d}$$

with $(t_k)_{k\in\mathbb{N}}$, $(\xi_k)_{k\in\mathbb{N}}$, $(\phi_k)_{k\in\mathbb{N}}$ and $(P_k)_{k\in\mathbb{N}}$ as given in the statement of Claim 1. The map $G_k(\cdot)$ is $C^0([t_{k-1},t_{k+1}]\times\mathbb{R}^n\times\mathbb{R}^{\bar{n}};\mathbb{R}^n)$ and locally Lipschitz on $(z,y)\in\mathbb{R}^n\times\mathbb{R}^{\bar{n}}$, thus, for any initial $z_k(t_0)\in\mathbb{R}^n$, $t_0 \in [t_{k-1},t_{k+1}]$, (74a) admits a unique solution in a neighborhood of $t_0$. Furthermore, $G_k:[t_{k-1},t_{k+1}]\times\mathbb{R}^n\times\mathbb{R}^{\bar{n}}\to\mathbb{R}^n$ is bounded, therefore the corresponding solution of (74) is defined for all $t \in [t_{k-1},t_{k+1}]$.

**Establishment of the convergence (5):** Let $Z:[t_0,\infty)\to\mathbb{R}^n$ as defined by (9), namely, $Z(t) := z_k(t)$, $t \in [t_k,t_{k+1})$, $k \in \mathbb{N}$, where for each $k \in \mathbb{N}$ the map $z_k(\cdot)$ is the solution of (74). Let $k_0 \in \mathbb{N}$ with $k_0 \geq 2$ such that for the initial state $x_0 \in \mathbb{R}^n$ of (11a) it holds:

$$k_0 \geq |x_0| \tag{75}$$

Consider now an arbitrary integer $k \geq k_0$. By invoking (74b) there exists a maximum time $\bar{t} \in (t_{k-1},t_{k+1}]$ for which $|z_k(t)| < \zeta_k$ for every $t \in [t_{k-1},\bar{t})$.

**Claim 2:** We claim that $\bar{t} = t_{k+1}$, namely:

$$|z_k(t)| < \zeta_k, \quad \forall t \in [t_{k-1},t_{k+1}) \tag{76}$$

**Proof of Claim 2:** Suppose on the contrary that there exists $\bar{t} \in (t_{k-1}, t_{k+1})$ with

$$|z_k(\bar{t})| = \zeta_k \text{ and } |z_k(t)| < \zeta_k, \ \forall t \in [t_{k-1}, \bar{t}) \tag{77}$$

By invoking (74a),(74c) and (77) it then follows that for the given $k$ above:

$$\dot{z}_k = F(t, z_k, y) + \phi_k(t) P_k^{-1}(t) H'(t)(y - H(t) z_k), \ \forall t \in [t_{k-1}, \bar{t}] \tag{78}$$

We now again recall Assumption A3, particularly, that both A1 and A2 hold with $R := k$ and notice that, due to definition (70) and (74b), the constant $\xi_k$ satisfies (19) with $R = k$ and $\bar{t}_0 := t_{k-1}$. Since $|x_0| \leq k = R$, it then follows, by taking into account (67),(68),(69),(78) and the result of Statement (ii) of Proposition 2.1, that estimation (21a) is valid, namely:

$$|e_k(t)| < \xi_k, \ \forall t \in [t_{k-1}, \bar{t}] \tag{79a}$$
$$\text{where } e_k(t) := x(t) - z_k(t), \ t \in [t_{k-1}, t_{k+1}] \tag{79b}$$

From (3),(74d),(75),(79) and by recalling that $k \geq k_0$ and $\bar{t} \leq t_{k+1}$, we estimate:

$$|z_k(\bar{t})| \leq |x(\bar{t})| + |e_k(\bar{t})| < \beta(\bar{t}, k) + \xi_k \leq \beta(t_{k+1}, k) + \xi_k = \zeta_k \tag{80}$$

which contradicts (77), therefore for the given $k$ (76) is fulfilled and the proof of Claim 2 is completed. ◀

It turns out from Claim 2 that $z_k(\cdot)$ satisfies (78) for every $t \in [t_{k-1}, t_{k+1}]$ which, by virtue of (21b) of Statement (ii) of Proposition 2.1 and (66a), implies:

$$|e_k(t)| \leq \max\{\xi_k \exp\{-c_1(t - (t_{k-1} + 1))\}, \sqrt{g(t)}\}, \ \forall t \in [t_{k-1} + 1, t_{k+1}] \tag{81}$$

By (71) and (81) we deduce:

$$|e_k(t)| \leq \frac{1}{k}, \ \forall t \in [t_k, t_{k+1}] \tag{82}$$

Finally, we show that the error $e(t) := x(t) - Z(t)$ between the solution of (11) and the map $Z(\cdot)$ as defined by (9), satisfies the estimation (5). Indeed, let an arbitrary $\varepsilon > 0$ and let $\bar{k} = \bar{k}(\varepsilon) \in \mathbb{N}$ with $\bar{k} \geq k_0$ such that $\frac{1}{\bar{k}} < \varepsilon$. Define $T = T(\varepsilon) := t_{\bar{k}}$. It follows from (66), that for every $t > T$ there exists an integer $k \geq \bar{k}$, in such a way that $t_k \leq t < t_{k+1}$ and therefore, since $k \geq k_0$, by (9),(79b) and (82) we have $|x(t) - Z(t)| = |x(t) - z_k(t)| = |e_k(t)| \leq \frac{1}{k} \leq \frac{1}{\bar{k}} < \varepsilon$ for every $t > T = T(\varepsilon)$, which implies (5) and the proof of Proposition 2.2 is completed. ■

## 3. Proof of Theorem 1.1

In this section we apply results of the previous section to prove our main result concerning the solvability of the SODP(ODP) for triangular systems (1).

**Proof of Theorem 1.1:** The establishment of both statements is based on the results of Propositions 2.1 and 2.2. In the analysis that follows, we may assume that, without any loss of generality, instead of assumption (10), it holds:

$$a_i(t, y) > 0, \quad \forall t \geq 0, \; y \in \mathbb{R}, \; i = 1, 2, \ldots, n-1 \tag{83}$$

For simplicity, we also may assume in the sequel that $m_i \geq 3$ for all $i \in \{1, 2, \ldots, n-1\}$. (The same analysis plus some elementary modifications can be used for the establishment of our statement for the case $m_i = 1$ for certain $i \in \{1, 2, \ldots, n-1\}$; the details are left to the reader). In order to prove the first statement, we show that Condition A3 is fulfilled for (1), namely, we prove that there exist an integer $\ell \in \mathbb{N}$, a map $A \in C^0(\mathbb{R}_{\geq 0} \times \mathbb{R}^\ell \times \mathbb{R}^n \times \mathbb{R}^n \times \mathbb{R}^k; \mathbb{R}^{n \times n})$ and constants $L > 1$, $c_1, c_2 > 0$ such that (14) holds and a function $g(\cdot)$ satisfying (13), in such a way that for each $R > 0$, both A1 and A2 hold for (1). Let $R > 0$, $\xi \geq 1$ and define:

$$F(t, x, y) := \begin{pmatrix} f_1(t, y) + a_1(t, y)x_2^{m_1} \\ f_2(t, y, x_2) + a_2(t, y)x_3^{m_2} \\ \vdots \\ f_n(t, y, x_2, \ldots, x_n) \end{pmatrix}, (t, x, y) \in \mathbb{R}_{\geq 0} \times \mathbb{R}^n \times \mathbb{R} \tag{84}$$

Let $\sigma_R := \sigma_{R,\xi} \in N \cap C^1([0, \infty); \mathbb{R})$ satisfying:

$$\sigma_R(t) \geq \sum_{i=2}^{n} \sum_{j=2}^{i} \max\left\{ \left| \frac{\partial f_i}{\partial x_j}(t, y, x_2, \ldots, x_i) \right| : |(y, x_2, \ldots, x_i)| \leq 2\beta(t, R) + \xi \right\}, \forall t \geq 0 \tag{85}$$

and consider the set-valued map $[0, \infty) \ni t \to Q_R(t) := Q_{R,\xi}(t) \subset \mathbb{R}^\ell$, $\ell := \frac{n(n+1)}{2}$ defined as

$$Q_R(t) := \{q = (q_{1,1}; q_{2,1}, q_{2,2}; \ldots; q_{n,1}, q_{n,2}, \ldots, q_{n,n}) \in \mathbb{R}^\ell : |q| \leq \sigma_R(t)\} \tag{86}$$

that obviously satisfies the CP. Also, let $Y_R(\cdot)$ as given by (12) with

$$H := (\underbrace{1, 0, \ldots, 0}_{n}) \tag{87}$$

From (84)-(87) and use of the mean-value theorem it follows that for every $t \geq 0$, $y \in Y_R(t)$ and $x, z \in \mathbb{R}^n$ with $|x| \leq \beta(t, R)$ and $|x - z| \leq \xi$ we have:

$$F(t, x, y) - F(t, z, y) = A(t, q, x, x-z, y)(x-z), \text{ for some } q \in Q_R(t) \text{ with } q_{i,1} = 0, i = 1, 2, \ldots, n;$$
(88a)

$$A(t, q, x, e, y) := \begin{pmatrix} q_{1,1} & a_1(t, y)b_1(x_2, e_2) & 0 & \cdots & 0 \\ q_{2,1} & q_{2,2} & a_2(t, y)b_2(x_3, e_3) & \ddots & \vdots \\ \vdots & \vdots & \vdots & \ddots & 0 \\ q_{n-1,1} & q_{n-1,2} & q_{n-1,3} & & a_{n-1}(t, y)b_{n-1}(x_n, e_n) \\ q_{n,1} & q_{n,2} & q_{n,3} & \cdots & q_{n,n} \end{pmatrix}$$
(88b)

$$b_i(\alpha,\beta) := \alpha^{m_i-1} + \alpha^{m_i-2}(\alpha-\beta) + \ldots + \alpha(\alpha-\beta)^{m_i-2} + (\alpha-\beta)^{m_i-1}, \quad (\alpha,\beta) \in \mathbb{R}^2, \quad i=1,2,\ldots,n-1 \tag{88c}$$

thus A1 is satisfied. Also, notice that, due to (83) and the fact that the map $Y_R$ as given by (12) satisfies the CP, there exists a function $w_R \in C^1(\mathbb{R}_{\geq 0}; \mathbb{R})$ such that

$$a_i(t,y) \geq w_R(t) > 0, \quad \forall t \geq 0, \quad y \in Y_R(t), \quad i=1,2,\ldots,n-1 \tag{89}$$

In the rest procedure we need the following elementary property concerning the mappings $b_i(\cdot,\cdot)$ as defined by (88c): For each $i=1,2,\ldots,n-1$ there exist a map $\bar{b}_i \in C(\mathbb{R}^2; \mathbb{R})$ and a constant $\vartheta_i > 0$ such that

$$b_i(\alpha,\beta) = \bar{b}_i(\alpha,\beta) + \vartheta_i \beta^{m_i-1}, \quad \forall (\alpha,\beta) \in \mathbb{R}^2; \tag{90a}$$

$$\bar{b}_i(\alpha,\beta) \geq 0, \quad \forall (\alpha,\beta) \in \mathbb{R}^2 \tag{90b}$$

Indeed, by taking into account that $m_i - 1 = 2\gamma_i$ for certain $\gamma_i \in \mathbb{N}$, it follows that $b_i(\alpha,\beta) > 0$ for every nonzero $(\alpha,\beta) \in \mathbb{R}^2$ and thus $\vartheta_i := \min\{b_i(\alpha,\beta) : (\alpha,\beta) \in \mathbb{R}^2 \text{ with } |(\alpha,\beta)|=1\} > 0$. Moreover, it holds $b_i(\lambda\alpha, \lambda\beta) = \lambda^{2\gamma_i} b_i((\alpha,\beta))$ for every $(\alpha,\beta) \in \mathbb{R}^2$ and constant $\lambda$. Therefore, if we define $\bar{b}_i(\alpha,\beta) := b_i(\alpha,\beta) - \vartheta_i \beta^{2\gamma_i}$, then for every nonzero $(\alpha,\beta) \in \mathbb{R}^2$ we have:

$$\bar{b}_i(\alpha,\beta) = b_i\left(\left(\frac{\alpha}{|(\alpha,\beta)|}, \frac{\beta}{|(\alpha,\beta)|}\right)\right)|(\alpha,\beta)| - \vartheta_i \beta^{2\gamma_i} = |(\alpha,\beta)|^{2\gamma_i} b_i\left(\frac{\alpha}{|(\alpha,\beta)|}, \frac{\beta}{|(\alpha,\beta)|}\right) - \vartheta_i \beta^{2\gamma_i}$$
$$\geq \vartheta_i((\alpha^2+\beta^2)^{\gamma_i} - \beta^{2\gamma_i}) \geq 0$$

hence, both (90a,b) are established.

Next, we establish A2, particularly, we prove that there exist constants $L>1$, $c_1, c_2 > 0$ such that (14) holds and a function $g(\cdot)$ satisfying (13), in such a way that for every $R>0$, $\xi \geq 1$ and $t_0 \geq 0$, a time-varying symmetric matrix $P_R := P_{R,\xi,t_0} \in C^1([t_0,\infty); \mathbb{R}^{n \times n})$ and a function $d_R := d_{R,\xi,t_0} \in C^0([t_0,\infty); \mathbb{R})$ can be found satisfying all conditions (17a,b,c) with $H$, $A(\cdot,\cdot,\cdot,\cdot,\cdot)$, $Y_R(\cdot)$ and $Q_R(\cdot)$, as given by (87),(88b),(12) and (86), respectively and with $\beta(\cdot,\cdot)$ as given in (3) for system (1). In order to establish A2, we proceed by induction as follows. Pick $L>1$, $c_1 := 1$, $c_2 := n$ (namely, $c_2$ is equal to the dimension of the state space) and let $g(\cdot)$ be a $C^1$ function satisfying (13) (for example we may take $g(t) := \frac{1}{2}e^{-t}$, $t \geq 0$). Also, let $R>0$ and for $k=2,3,\ldots,n$ define:

$$H_k := (\underbrace{1,0,\ldots,0}_{k}), \quad e := (e_{n-k+1}; \hat{e}')' \in \mathbb{R} \times \mathbb{R}^{k-1}, \quad \hat{e} := (e_{n-k+2},\ldots,e_n)' \in \mathbb{R}^{k-1}; \tag{91a}$$

$$A_k(t,q,x,e,y) := \begin{pmatrix} q_{n-k+1,n-k+1} & a_{n-k+1}(t,y)b_{n-k+1}(x_{n-k+2},e_{n-k+2}) & 0 & \cdots & 0 \\ q_{n-k+2,n-k+1} & & & & \\ \vdots & & \boxed{A_{k-1}(t,q,x,\hat{e},y)} & & \\ q_{n,n-k+1} & & & & \end{pmatrix} \tag{91b}$$

where

$$A_1(t,q,x,e_n,y) := q_{n,n} \tag{91c}$$

**Claim 1 (Induction Hypothesis):** Let $k \in \mathbb{N}$ with $2 \le k \le n$. Then for $L$, $R$ and $g(\cdot)$ as above and for every $\xi \ge 1$ and $t_0 \ge 0$, there exist a time-varying symmetric matrix $P_{R,k} := P_{R,\xi,t_0,k} \in C^1([t_0,\infty); \mathbb{R}^{k \times k})$ and a mapping $d_{R,k} := d_{R,\xi,t_0,k} \in C^0([t_0,\infty); \mathbb{R})$, in such a way that the following hold:

$$P_{R,k}(t) > I_{k \times k}, \ \forall t \ge t_0 \ ; \ |P_{R,k}(t_0)| \le L; \tag{92a}$$

$$d_{R,k}(t) > n-k+1, \ \forall t \ge t_0+1; \ \int_{t_1}^{t_2} d_{R,k}(s)ds > -k, \ \forall t_2 \ge t_1, \ t_1,t_2 \in [t_0, t_0+1]; \tag{92b}$$

$$e'P_{R,k}(t)A_k(t,q,x,e,y)e + \tfrac{1}{2}e'\dot{P}_{R,k}(t)e \le -d_{R,k}(t)e'P_{R,k}(t)e,$$
$$\forall t \ge t_0, \ q \in Q_R(t), \ x \in \mathbb{R}^n, \ e \in \ker H_k, \ y \in Y_R(t): \tag{92c}$$
$$|x| \le \beta(t,R), \ |e| \le \xi, \ e'P_{R,k}(t)e \ge g(t)$$

with $H_k$, $A_k(\cdot,\cdot,\cdot,\cdot,\cdot)$, $Y_R(\cdot)$ and $Q_R(\cdot)$ as given in (91a),(91b),(12) and (86), respectively.

Notice that, according to (91), the mappings $H_n$ and $A_n(\cdot,\cdot,\cdot,\cdot,\cdot)$ coincide with $H$ and $A(\cdot,\cdot,\cdot,\cdot,\cdot)$ as given by (87) and (88b), respectively and A2 is a consequence of Claim 1 with $H := H_n$ and $A(\cdot,\cdot,\cdot,\cdot,\cdot) := A_n(\cdot,\cdot,\cdot,\cdot,\cdot)$ and with $d_R := d_{R,n}$ and $P_R := P_{R,n}$ as given in (92a,b). Indeed, relations (17a),(17c) follow directly from (92a),(92c) and both inequalities of (17b) follow from (92b) with $k = n$ and by taking into account that $c_1 = 1$ and $c_2 = n$.

**Proof of Claim 1 for** $k := 2$: For the case $k = 2$ let

$$H_2 := (1,0), \ e := (e_{n-1}, e_n)' \in \mathbb{R}^2 \tag{93a}$$

$$A_2(t,q,x,e,y) := \begin{pmatrix} q_{n-1,n-1} & a_{n-1}(t,y)b_{n-1}(x_n,e_n) \\ q_{n,n-1} & q_{n,n} \end{pmatrix} \tag{93b}$$

For $L$, $R$ and $g(\cdot)$ as above and for every $\xi \ge 1$ and $t_0 \ge 0$, we establish existence of a time-varying symmetric matrix $P_{R,2} := P_{R,\xi,t_0,2} \in C^1([t_0,\infty); \mathbb{R}^{2 \times 2})$ and a mapping $d_{R,2} := d_{R,\xi,t_0,2} \in C^0([t_0,\infty); \mathbb{R})$ in such a way that

$$P_{R,2}(t) > I_{2 \times 2}, \ \forall t \ge t_0 \ ; \ |P_{R,2}(t_0)| \le L; \tag{94a}$$

$$d_{R,2}(t) > n-1, \forall t \ge t_0+1; \int_{t_1}^{t_2} d_{R,2}(s)ds > -2, \forall t_2 \ge t_1, \ t_1,t_2 \in [t_0, t_0+1]; \tag{94b}$$

$$e'P_{R,2}(t)A_2(t,q,x,e,y)e + \tfrac{1}{2}e'\dot{P}_{R,2}(t)e \le -d_{R,2}(t)e'P_{R,2}(t)e,$$
$$\forall t \ge t_0, \ q \in Q_R(t), \ x \in \mathbb{R}^n, \ e = (e_{n-1}, e_n)' \in \mathbb{R}^2, \ y \in Y_R(t): \tag{94c}$$
$$|x| \le \beta(t,R), \ e \in \ker H_2, \ |e| \le \xi, \ e'P_{R,2}(t)e \ge g(t)$$

with $H_2$, $A_2(\cdot,\cdot,\cdot,\cdot,\cdot)$, $Y_R(\cdot)$ and $Q_R(\cdot)$ as given in (93a),(93b),(12) and (86), respectively. Let $\xi \ge 1$, $t_0 \ge 0$ and define:

$$P_{R,2}(t) := \begin{pmatrix} p_{R,1}(t) & p_R(t) \\ p_R(t) & L \end{pmatrix}, \ t \ge t_0 \tag{95}$$

for certain $p_{R,1}, p_R \in C^1([t_0,\infty);\mathbb{R})$, yet to be determined. Notice that, due to (93a) and (95), we have:

$$\begin{aligned}&\{e \in \ker H_2 : |e| \leq \xi \text{ and } e'P_{R,2}(t)e \geq g(t)\}\\&= \{e = (0, e_n)' : \sqrt{g(t)/L} \leq |e_n| \leq \xi\}\end{aligned} \quad (96)$$

Then, by taking into account (93),(95) and (96), the desired (94c) is written:

$$(0, e_n)\begin{pmatrix} p_{R,1}(t) & p_R(t) \\ p_R(t) & L \end{pmatrix}\begin{pmatrix} q_{n-1,n-1} & a_{n-1}(t,y)b_{n-1}(x_n, e_n) \\ q_{n,n-1} & q_{n,n} \end{pmatrix}\begin{pmatrix} 0 \\ e_n \end{pmatrix} + \frac{1}{2}(0, e_n)\overline{\begin{pmatrix} p_{R,1}(t) & p_R(t) \\ p_R(t) & L \end{pmatrix}}\begin{pmatrix} 0 \\ e_n \end{pmatrix}$$

$$\leq -d_{R,2}(t)(0, e_n)\begin{pmatrix} p_{R,1}(t) & p_R(t) \\ p_R(t) & L \end{pmatrix}\begin{pmatrix} 0 \\ e_n \end{pmatrix}$$

$$\forall t \geq t_0,\, q \in Q_R(t),\, x \in \mathbb{R}^n,\, e = (0, e_n)' \in \mathbb{R}^2,\, y \in Y_R(t):$$
$$|x| \leq \beta(t, R),\, \sqrt{g(t)/L} \leq |e_n| \leq \xi$$

or equivalently:

$$p_R(t)a_{n-1}(t,y)b_{n-1}(x_n,e_n) + Lq_{n,n} \leq -Ld_{R,2}(t),$$
$$\forall t \geq t_0,\, q \in Q_R(t),\, x \in \mathbb{R}^n,\, e = (0, e_n)' \in \mathbb{R}^2,\, y \in Y_R(t): \quad (97)$$
$$|x| \leq \beta(t, R),\, \sqrt{g(t)/L} \leq |e_n| \leq \xi$$

We impose the additional requirement for the candidate $p_R(\cdot)$:

$$p_R(t) \leq 0,\, \forall t \geq t_0\,;\, p_R(t_0) = 0 \quad (98)$$

By taking into account (86),(89),(90),(98) and the equivalence between (94c) and (97), it follows that, in order to prove (94c), it suffices to determine $p_{R,1}, p_R \in C^1([t_0,\infty);\mathbb{R})$ and $d_{R,2} \in C^0([t_0,\infty);\mathbb{R})$ in such a way that (94a,b) and (98) are fulfilled, and further:

$$p_R(t)w_R(t)(\overline{b}_{n-1}(x_n, e_n) + \vartheta_{n-1}e_n^{m_{n-1}-1}) + L\sigma_R(t) \leq -Ld_{R,2}(t),$$
$$\forall t \geq t_0,\, x \in \mathbb{R}^n,\, e_n \in \mathbb{R}: \quad (99)$$
$$|x| \leq \beta(t, R),\, \sqrt{g(t)/L} \leq |e_n| \leq \xi$$

Then, by taking into account (89),(90b),(98) and (99), it suffices to prove:

$$p_R(t)w_R(t)\vartheta_{n-1}\frac{g^{\frac{m_{n-1}-1}{2}}(t)}{L^{\frac{m_{n-1}-1}{2}}} + L\sigma_R(t) \leq -Ld_{R,2}(t),\, \forall t \geq t_0 \quad (100)$$

for certain $p_{R,1}, p_R \in C^1([t_0,\infty);\mathbb{R})$ and $d_{R,2} \in C^0([t_0,\infty);\mathbb{R})$ satisfying (94a,b) and (98).

**Construction of $p_R$ and $d_{R,2}$:** Let

$$M_2 := \max\{\sigma_R(t) : t \in [t_0, t_0 + \tfrac{1}{2}]\} \quad (101a)$$

$$\tau_2 := \min\left\{\frac{1}{M_2}, 1\right\} \quad (101b)$$

and define $\theta := \theta_{R,\xi,t_0} \in C^1([t_0,\infty);\mathbb{R})$, $p_R \in C^1([t_0,\infty);\mathbb{R})$ and $d_{R,2} \in C^0([t_0,\infty);\mathbb{R})$ as follows:

$$\theta(t) \begin{cases} := 0, & t = t_0 \\ \in [0,1], & t \in [t_0, t_0 + \frac{\tau_2}{2}] \\ := 1, & t \in [t_0 + \frac{\tau_2}{2}, \infty) \end{cases} \quad (102)$$

$$p_R(t) := -\theta(t) \frac{L^{\frac{m_{n-1}+1}{2}} (n + \sigma_R(t))}{\vartheta_{n-1} w_R(t) g^{\frac{m_{n-1}-1}{2}}(t)}, \quad t \geq t_0 \quad (103)$$

$$d_{R,2}(t) \begin{cases} := -M_2, & t \in [t_0, t_0 + \frac{\tau_2}{2}] \\ \in [-M_2, n], & t \in [t_0 + \frac{\tau_2}{2}, t_0 + \tau_2] \\ := n, & t \in [t_0 + \tau_2, \infty) \end{cases} \quad (104)$$

We show that (94b),(98) and (100) are fulfilled, with $p_R(\cdot)$ and $d_{R,2}(\cdot)$ as given by (103) and (104), respectively. Indeed, (98) follows directly by invoking (89),(102),(103) and recalling that $\vartheta_{n-1} > 0$. The first inequality of (94b) is also a direct consequence of (101b) and (104). Also, by (101b) and (104) we have $\int_{t_1}^{t_2} d_{R,2}(s) ds \geq -\tau_2 M_2 \geq -1$ for every $t_2 \geq t_1$, $t_1, t_2 \in [t_0, t_0 + 1]$, hence, the second inequality of (94b) is fulfilled. We next show that (100) holds as well, with $p_R(\cdot)$ and $d_{R,2}(\cdot)$ as previously defined. We consider two cases:

**Case 1:** $t \in [t_0, t_0 + \frac{\tau_2}{2}]$. Then, by invoking (89),(98),(101),(104) and recalling that $\vartheta_{n-1} > 0$ it follows that

$$p_R(t) w_R(t) \vartheta_{n-1} \frac{g^{\frac{m_{n-1}-1}{2}}(t)}{L^{\frac{m_{n-1}-1}{2}}} + L\sigma_R(t) \leq L\sigma_R(t) \leq LM_2 = -Ld_{R,2}(t)$$

hence, (100) is satisfied for every $t \in [t_0, t_0 + \frac{\tau_2}{2}]$.

**Case 2:** $t \in [t_0 + \frac{\tau_2}{2}, \infty)$. Then from (102),(103) and (104) we evaluate:

$$p_R(t) w_R(t) \vartheta_{n-1} \frac{g^{\frac{m_{n-1}-1}{2}}(t)}{L^{\frac{m_{n-1}-1}{2}}} + L\sigma_R(t) = -Ln \leq -Ld_{R,2}(t)$$

namely, (100) again holds for all $t \in [t_0 + \frac{\tau_2}{2}, \infty)$.

We therefore conclude that (100) is fulfilled for all $t \geq t_0$.

**Construction of** $p_{R,1}$: It remains to determine a function $p_{R,1} \in C^1([t_0,\infty);\mathbb{R})$ such that both requirements of (94a) are fulfilled, with $P_{R,2}(\cdot)$ as given by (95) and (103). By recalling that $L > 1$, we may define:

$$p_{R,1}(t) := \frac{p_R^2(t)}{L-1} + L, \quad t \geq t_0 \quad (105a)$$

with $p_R(\cdot)$ as given by (103). From (95) and (105a) we deduce:

$$\det(P_{R,2}(t) - I_{2\times 2}) = \det\begin{pmatrix} p_{R,1}(t) - 1 & p_R(t) \\ p_R(t) & L-1 \end{pmatrix} = (L-1)^2 > 0, \quad \forall t \geq t_0 \quad (105b)$$

From (105b) and our assumption $L>1$ it follows that $P_{R,2}(t) - I_{2\times 2} > 0$ for every $t \geq t_0$, hence, the first inequality of (94a) is fulfilled. Finally, the second inequality of (94a) follows directly from (95),(98) and (105a). This completes the proof of Claim 1 for $k=2$.

**Proof of Claim 1 (general step of induction procedure):** Assume now that the statement of Claim 1 is fulfilled for certain integer $k$ with $2 \leq k < n$; particularly, assume that for $L$, $R$ and $g(\cdot)$ as above and for all $\xi \geq 1$ and $t_0 \geq 0$ there exist a time-varying symmetric matrix $P_{R,k} := P_{R,\xi,t_0,k} \in C^1([t_0,\infty);\mathbb{R}^{k\times k})$ and a mapping $d_{R,k} := d_{R,\xi,t_0,k} \in C^0([t_0,\infty);\mathbb{R})$, such that (92a,b,c) are fulfilled, with $H_k$, $A_k(\cdot,\cdot,\cdot,\cdot,\cdot)$, $Y_R(\cdot)$ and $Q_R(\cdot)$, as given in (91a),(91b),(12) and (86), respectively. We prove that Claim 1 also holds for $k := k+1$. Consider the pair $(H,A)$ as given in (22) with

$$H(t) := H_k, \; A(t,q,x,e,y) := A_k(t,q,x,e,y), \; \ell = \frac{n(n+1)}{2}, \; m := k, \; n := n \text{ and } \bar{n} := 1$$

where $H_k$ and $A_k$ are defined by (91a) and (91b), respectively. Then, due to the first inequality of (92a) and (92c), Hypothesis H1 of the previous section is fulfilled, with $g(\cdot)$ as above, $Y_R(\cdot)$, $Q_R(\cdot)$ and $\beta_R(\cdot) := \beta(\cdot,R)$ as given in (12),(86) and (3), respectively, and with $d_R(\cdot) := d_{R,k}(\cdot)$ and $P_R(\cdot) := P_{R,k}(\cdot)$ as given in (92a,b). Then, if we consider the function $\bar{d}_{R,k} := \bar{d}_{R,\xi,t_0,k}$ defined as:

$$\bar{d}_{R,k}(t) := d_{R,k}(t) - \tfrac{1}{2}, \; t \geq t_0 \tag{106}$$

we have $\bar{d}_{R,k}(t) < d_{R,k}(t)$ for all $t \geq t_0$, consequently, all requirements of Lemma 2.1 are fulfilled and therefore, there exists a function $\phi_{R,k} := \phi_{R,\xi,t_0,k} \in C^1([t_0,\infty);\mathbb{R}_{>0})$ such that

$$e'P_{R,k}(t)A_k(t,q,x,e,y)e + \tfrac{1}{2}e'\dot{P}_{R,k}(t)e \leq \phi_{R,k}(t)|H_k e|^2 - \bar{d}_{R,k}(t)e'P_{R,k}(t)e,$$
$$\forall t \geq t_0, q \in Q_R(t), x \in \mathbb{R}^n, e \in \mathbb{R}^k, y \in Y_R(t): \tag{107}$$
$$|x| \leq \beta(t,R), |e| \leq \xi, e'P_{R,k}(t)e \geq g(t)$$

Notice that, due to (92b) and (106), $\bar{d}_{R,k}(\cdot)$ satisfies:

$$\bar{d}_{R,k}(t) > n-k+\tfrac{1}{2}, \; \forall t \geq t_0+1; \; \int_{t_1}^{t_2} \bar{d}_{R,k}(s)ds > -(k+\tfrac{1}{2}), \; \forall t_2 \geq t_1, \; t_1,t_2 \in [t_0,t_0+1] \tag{108}$$

By exploiting (107) and (108), we are in a position to establish that Claim 1 is fulfilled for $k = k+1$. Specifically, for the same $L$, $R$ and $g(\cdot)$ as above and for any $\xi \geq 1$ and $t_0 \geq 0$, we show that there exist a time-varying symmetric matrix $P_{R,k+1} \in C^1([t_0,\infty);\mathbb{R}^{(k+1)\times(k+1)})$ and a map $d_{R,k+1} \in C^0([t_0,\infty);\mathbb{R})$, such that (92a,b) are fulfilled with $k = k+1$ and further:

$$e'P_{R,k+1}(t)A_{k+1}(t,q,x,e,y)e + \tfrac{1}{2}e'\dot{P}_{R,k+1}(t)e \leq -d_{R,k+1}(t)e'P_{R,k+1}(t)e,$$
$$\forall t \geq t_0, q \in Q_R(t), x \in \mathbb{R}^n, e \in \ker H_{k+1}, y \in Y_R(t): \tag{109}$$
$$|x| \leq \beta(t,R), |e| \leq \xi, e'P_{R,k+1}(t)e \geq g(t)$$

where

$$H_{k+1} := \underbrace{(1,0,...,0)}_{k+1}, \; e := (e_{n-k}; \hat{e}')' \in \mathbb{R} \times \mathbb{R}^k, \hat{e} := (e_{n-k+1},...,e_n)' \in \mathbb{R}^k; \tag{110a}$$

$$A_{k+1}(t,q,x,e,y) := \begin{pmatrix} q_{n-k,n-k} & a_{n-k}(t,y)b_{n-k}(x_{n-k+1},e_{n-k+1}) & 0 & \cdots & 0 \\ q_{n-k+1,n-k} & & & & \\ \vdots & & \boxed{A_k(t,q,x,\hat{e},y)} & & \\ q_{n,n-k} & & & & \end{pmatrix}; \quad (110b)$$

$$P_{R,k+1}(t) := \begin{pmatrix} p_{R,1}(t) & p_R(t) & 0 & \cdots & 0 \\ p_R(t) & & & & \\ 0 & & \boxed{P_{R,k}(t)} & & \\ \vdots & & & & \\ 0 & & & & \end{pmatrix} \quad (110c)$$

and where $Y_R(\cdot)$, $Q_R(\cdot)$ are given in (12) and (86), respectively. For given $\xi \geq 1$ and $t_0 \geq 0$ we determine $p_{R,1}, p_R \in C^1([t_0,\infty);\mathbb{R})$ and $d_{R,k+1} \in C^0([t_0,\infty);\mathbb{R})$ such that (92a,b) are fulfilled with $k = k+1$, and further (109) holds, with $H_{k+1}$, $A_{k+1}(\cdot,\cdot,\cdot,\cdot,\cdot)$ and $P_{R,k+1}(\cdot)$ as given by (110), namely:

$$(0, e_{n-k+1}, \ldots, e_n) \begin{pmatrix} p_{R,1}(t) & p_R(t) & 0 & \cdots & 0 \\ p_R(t) & & & & \\ 0 & & \boxed{P_{R,k}(t)} & & \\ \vdots & & & & \\ 0 & & & & \end{pmatrix} \begin{pmatrix} q_{n-k,n-k} & a_{n-k}(t,y)b_{n-k}(x_{n-k+1},e_{n-k+1}) & 0 & \cdots & 0 \\ q_{n-k+1,n-k} & & & & \\ \vdots & & \boxed{A_k(t,q,x,\hat{e},y)} & & \\ q_{n,n-k} & & & & \end{pmatrix} \begin{pmatrix} 0 \\ e_{n-k+1} \\ \vdots \\ e_n \end{pmatrix}$$

$$+ \tfrac{1}{2}(0, e_{n-k+1}, \ldots, e_n) \overbrace{\begin{pmatrix} p_{R,1}(t) & p_R(t) & 0 & \cdots & 0 \\ p_R(t) & & & & \\ 0 & & \boxed{P_{R,k}(t)} & & \\ \vdots & & & & \\ 0 & & & & \end{pmatrix}}^{\displaystyle \cdot} \begin{pmatrix} 0 \\ e_{n-k+1} \\ \vdots \\ e_n \end{pmatrix}$$

$$\leq -d_{R,k+1}(t)(0, e_{n-k+1}, \ldots, e_n) \begin{pmatrix} p_{R,1}(t) & p_R(t) & 0 & \cdots & 0 \\ p_R(t) & & & & \\ 0 & & \boxed{P_{R,k}(t)} & & \\ \vdots & & & & \\ 0 & & & & \end{pmatrix} \begin{pmatrix} 0 \\ e_{n-k+1} \\ \vdots \\ e_n \end{pmatrix},$$

$$\forall t \geq t_0, q \in Q_R(t), x \in \mathbb{R}^n, e := (e_{n-k}; \hat{e}')' \in \mathbb{R} \times \mathbb{R}^k, y \in Y_R(t):$$
$$|x| \leq \beta(t,R), e \in \ker H_{k+1}, |e| \leq \xi, e'P_{R,k+1}(t)e \geq g(t)$$

or, equivalently:

$$e_{n-k+1}^2 p_R(t)a_{n-k}(t,y)b_{n-k}(x_{n-k+1},e_{n-k+1}) + \hat{e}'P_{R,k}(t)A_k(t,q,x,\hat{e},y)\hat{e} + \tfrac{1}{2}\hat{e}'\dot{P}_{R,k}(t)\hat{e}$$
$$\leq -d_{R,k+1}\hat{e}'P_{R,k}(t)\hat{e}, \quad (111)$$
$$\forall t \geq t_0, q \in Q_R(t), x \in \mathbb{R}^n, e := (e_{n-k}; \hat{e}')' \in \mathbb{R} \times \mathbb{R}^k, y \in Y_R(t):$$
$$|x| \leq \beta(t,R), e \in \ker H_{k+1}, |e| \leq \xi, e'P_{R,k+1}(t)e \geq g(t)$$

Notice that, according to (110a,c), we have $e'P_{R,k+1}(t)e = \hat{e}'P_{R,k}(t)\hat{e}$ for every $e = (e_{n-k}; \hat{e}')' \in \ker H_{k+1}$, thus, by taking into account (91a) and (107), it follows that in order to prove (111), it suffices to show that

$$e^2_{n-k+1}(p_R(t)a_{n-k}(t,y)b_{n-k}(x_{n-k+1},e_{n-k+1}) + \phi_{R,k}(t)) \le (\bar{d}_{R,k}(t) - d_{R,k+1}(t))\hat{e}'P_{R,k}(t)\hat{e},$$
$$\forall t \ge t_0,\, x \in \mathbb{R}^n,\, \hat{e} \in \mathbb{R}^k,\, y \in Y_R(t): \qquad (112)$$
$$|x| \le \beta(t,R),\, |\hat{e}| \le \xi,\, \hat{e}'P_{R,k}(t)\hat{e} \ge g(t)$$

We impose the following additional requirements for the candidate functions $p_R(\cdot)$ and $d_{R,k+1}(\cdot)$:

$$p_R(t) \le 0,\, \forall t \ge t_0;\, p_R(t_0) = 0; \qquad (113\mathrm{a})$$
$$d_{R,k+1}(t) \le \bar{d}_{R,k}(t),\, \forall t \ge t_0 \qquad (113\mathrm{b})$$

By taking into account (89),(90),(113) and the fact that desired inequality in (112) should be valid for those $\hat{e} \in \mathbb{R}^k$ for which $|\hat{e}| \le \xi$ and $\hat{e}'P_{R,k}(t)\hat{e} \ge g(t)$, it follows that in order to show that (112) together with (92a,b) are valid with $k = k+1$, it suffices to show that

$$e^2_{n-k+1}(p_R(t)w_R(t)\vartheta_{n-k}e^{m_{n-k}-1}_{n-k+1} + \phi_{R,k}(t)) \le (\bar{d}_{R,k}(t) - d_{R,k+1}(t))g(t),$$
$$\forall t \ge t_0,\, |e_{n-k+1}| \le \xi \qquad (114)$$

for appropriate $p_{R,1}, p_R \in C^1([t_0,\infty);\mathbb{R})$ and $d_{R,k+1} \in C^0([t_0,\infty);\mathbb{R})$ in such a way that (92a,b) hold with $k = k+1$, and in addition $p_R(\cdot)$ and $d_{R,k+1}(\cdot)$ satisfy (113). Next, we proceed to the explicit construction of these functions.

**Construction of the map $p_R$:** Consider the integer $\gamma := \dfrac{m_{n-k}+1}{2}$ and define:

$$p_R(t) := -\frac{\zeta(t)\phi^\gamma_{R,k}(t)}{\vartheta_{n-k}w_R(t)},\, t \ge t_0 \qquad (115)$$

for certain $\zeta(\cdot) := \zeta_{R,\xi,t_0}(\cdot)$ satisfying:

$$\zeta(t) \ge 0,\, \forall t \ge t_0;\, \zeta(t_0) = 0 \qquad (116)$$

yet to be determined. It turns out that, in order to show (114) with $p_R(\cdot)$ and $d_{R,k+1}(\cdot)$ as above, the analysis is reduced to the construction of a pair of functions $\zeta \in C^1([t_0,\infty);\mathbb{R}_{\ge 0})$ and $d_{R,k+1} \in C^0([t_0,\infty);\mathbb{R})$ satisfying (113b),(116) and (92b) with $k = k+1$, and further:

$$-\zeta(t)s^\gamma + s \le (\bar{d}_{R,k}(t) - d_{R,k+1}(t))g(t),\, \forall t \ge t_0,\, s \in [0,\phi_{R,k}(t)\xi^2] \qquad (117)$$

**Construction of the mappings $\zeta$ and $d_{R,k+1}$:** Let

$$M_{k+1} := \left|\max\left\{-\bar{d}_{R,k}(t) + \frac{1}{g(t)}s + \frac{1}{4} : t \in [t_0, t_0 + \tfrac{1}{2}],\, s \in [0,\phi_{R,k}(t)\xi^2]\right\}\right| \qquad (118\mathrm{a})$$

$$\tau_{k+1} := \min\left\{\frac{1}{4M_{k+1}}, \frac{1}{2}\right\} \qquad (118\mathrm{b})$$

and consider a function $\theta := \theta_{R,\xi,t_0} \in C^1([t_0,\infty);\mathbb{R})$ defined as:

$$\theta(t) \begin{cases} := 0, & t = t_0 \\ \in [0,1], & t \in [t_0, t_0 + \frac{\tau_{k+1}}{2}] \\ := 1, & t \in [t_0 + \frac{\tau_{k+1}}{2}, \infty) \end{cases} \quad (119)$$

Notice that, due to (118), it holds:

$$\bar{d}_{R,k}(t) - \tfrac{1}{4} \geq -M_{k+1}, \quad \forall t \in [t_0, t_0 + \tau_{k+1}] \quad (120)$$

By taking into account (120), a function $d_{R,k+1} \in C^0([t_0, \infty); \mathbb{R})$ can be constructed satisfying:

$$d_{R,k+1}(t) \begin{cases} := -M_{k+1}, & t \in [t_0, t_0 + \frac{\tau_{k+1}}{2}] \\ \in [-M_{k+1}, \bar{d}_{R,k}(t) - \tfrac{1}{4}], & t \in [t_0 + \frac{\tau_{k+1}}{2}, t_0 + \tau_{k+1}] \\ := \bar{d}_{R,k}(t) - \tfrac{1}{4}, & t \in [t_0 + \tau_{k+1}, \infty) \end{cases} \quad (121)$$

Finally, define:

$$\zeta(t) := \frac{\theta(t)}{\gamma(\tfrac{1}{4} g(t))^{\gamma-1}}, \quad t \geq t_0 \quad (122)$$

that obviously is of class $C^1([t_0, \infty); \mathbb{R})$ and due to (119), satisfies (116). We next show that the functions $\zeta$ and $d_{R,k+1}$ satisfy the desired property (92b) for $k = k+1$ and the inequality (117).

**Proof of (92b) for $k = k+1$:** First, by virtue of (108),(118b) and (121) we can easily deduce that $d_{R,k+1}(t) > n - k + \tfrac{1}{4} > n - (k+1) + 1$ for all $t \geq t_0 + 1$, thus, the first inequality of (92b) holds for $k = k+1$. Also, by (118b) and (121) we have:

$$t_1, t_2 \in [t_0, t_0 + \tau_{k+1}], t_2 \geq t_1 \Rightarrow \int_{t_1}^{t_2} d_{R,k+1}(s) ds \geq -\int_{t_1}^{t_2} M_{k+1} ds \geq -\int_{t_0}^{t_0 + \tau_{k+1}} M_{k+1} ds \geq -\tfrac{1}{4} \quad (123a)$$

and by using (108),(118b) and (121) we get:

$$t_1, t_2 \in [t_0 + \tau_{k+1}, t_0 + 1], t_2 \geq t_1 \Rightarrow \int_{t_1}^{t_2} d_{R,k+1}(s) ds = \int_{t_1}^{t_2} \bar{d}_{R,k}(s) ds + \int_{t_1}^{t_2} (d_{R,k+1}(s) - \bar{d}_{R,k}(s)) ds$$
$$> -(k + \tfrac{1}{2}) - \int_{t_1}^{t_2} \tfrac{1}{4} ds \geq -(k + \tfrac{1}{2}) - \tfrac{1}{4} = -(k + \tfrac{3}{4})$$
$$(123b)$$

From (123) it follows that the second inequality of (92b) is fulfilled as well for $k = k+1$, namely, it holds:

$$\int_{t_1}^{t_2} d_{R,k+1}(s) ds > -(k+1), \quad \forall t_2 \geq t_1, \, t_1, t_2 \in [t_0, t_0 + 1]$$

**Proof of (117):** We consider two cases:

**Case 1:** $t \in [t_0, t_0 + \frac{\tau_{k+1}}{2}]$. First, notice that (118a) asserts that $M_{k+1} \geq -\bar{d}_{R,k}(t) + \dfrac{\phi_{R,k}(t)\xi^2}{g(t)}$ for every $t \in [t_0, t_0 + \tfrac{1}{2}]$, which in conjunction with (118b) and (121) imply:

$$d_{R,k+1}(t) \leq \bar{d}_{R,k}(t) - \frac{\phi_{R,k}(t)\xi^2}{g(t)}, \quad \forall t \in [t_0, t_0 + \frac{\tau_{k+1}}{2}] \qquad (124)$$

It follows by taking into account (116) and (124) that $-\zeta(t)s^\gamma + s - (\bar{d}_{R,k}(t) - d_{R,k+1}(t))g(t) \leq s - (\bar{d}_{R,k}(t) - d_{R,k+1}(t))g(t) \leq s - \phi_{R,k}(t)\xi^2 \leq 0$ for all $s \in [0, \phi_{R,k}(t)\xi^2]$, which implies (117).

**Case 2:** $t \in [t_0 + \frac{\tau_{k+1}}{2}, \infty)$. Notice that, due to (118b) and (121), we have $\bar{d}_{R,k}(t) - d_{R,k+1}(t) \geq \frac{1}{4}$ for every $t \in [t_0 + \frac{\tau_{k+1}}{2}, \infty)$, hence, by taking into account (119) and (122), it follows that

$$-\zeta(t)s^\gamma + s - (\bar{d}_{R,k}(t) - d_{R,k+1}(t))g(t) \leq -\zeta(t)s^\gamma + s - \tfrac{1}{4}g(t) =$$
$$-\frac{1}{\gamma(\tfrac{1}{4}g(t))^{\gamma-1}}s^\gamma + s - \tfrac{1}{4}g(t), \forall t \in [t_0 + \tfrac{\tau_{k+1}}{2}, \infty), s \geq 0 \qquad (125)$$

Also, it can be easily verified that for each fixed $t \geq 0$ the values of the function

$$[0, \infty) \ni s \to -\frac{s^\gamma}{\gamma(\tfrac{1}{4}g(t))^{\gamma-1}} + s - \tfrac{1}{4}g(t)$$

are strictly negative. The latter, in conjunction with (125) asserts that (117) is fulfilled for every $t \in [t_0 + \frac{\tau_{k+1}}{2}, \infty)$.

Both cases above guarantee that (117) holds for all $t \in [t_0, \infty)$ as required.

**Construction of the map** $p_{R,1}$: It remains to determine a function $p_{R,1} \in C^1([t_0, \infty); \mathbb{R})$ such that both inequalities of (92a) hold, with $P_{R,k+1}(\cdot)$ as given by (110c) and (115). Due to the first inequality of (92a), we may define:

$$p_{R,1}(t) := \frac{p_R^2(t) \det(P_{R,k-1}(t) - I_{(k-1)\times(k-1)})}{\det(P_{R,k}(t) - I_{k\times k})} + L, \quad t \geq t_0 \qquad (126a)$$

Also, since $L > 1$ and by invoking the first inequality of (92a) it follows from (110c) and (126a):

$$\det(P_{R,k+1}(t) - I_{(k+1)\times(k+1)}) = \det\begin{pmatrix} p_{R,1}(t)-1 & p_R(t) & 0 & \cdots & 0 \\ p_R(t) & & & & \\ 0 & & \boxed{P_{R,k}(t) - I_{k\times k}} & & \\ \vdots & & & & \\ 0 & & & & \end{pmatrix} = (L-1)\det(P_{R,k}(t) - I_{k\times k}) > 0, \forall t \geq t_0$$

(126b)

We deduce, from the first inequality of (92a), (126b) and the Sylvester criterion, that $P_{R,k+1}(t) > I_{(k+1)\times(k+1)}$ for every $t \geq t_0$, hence, the first inequality of (92a) is fulfilled for $k = k+1$. Finally, the second inequality of (92a) for $P_{R,k+1}(\cdot)$ is a direct consequence of (110c), (113a), (126a) and the induction hypothesis (92a); particularly, that $P_k(\cdot)$ satisfies $|P_k(t_0)| \leq L$.

We have shown that all requirements of Claim 1 hold, which, as was pointed out, establishes that for every $R > 0$ Hypothesis A2 is fulfilled. We therefore conclude that for system (1) Assumption A3 is satisfied hence, by invoking the result of Proposition 2.2, it follows that the SODP is solvable for (1) with respect to $\mathbb{R}^n$. The establishment of the second statement of Theorem 1.1, namely, for the case where it is known that the initial states of (1) belong to a ball $B_R$ of radius $R > 0$ centered at $0 \in \mathbb{R}^n$, follows directly from Proposition 2.1. For completeness we note that in that case, we may apply the same analysis, in order to establish that both A1 and A2 are fulfilled for the specific constant $R$, and then we may use the result of Proposition 2.1 to prove that the ODP is solvable for (1) with respect to $B_R$. ∎

## 4 Conclusion

For a class of triangular nonlinear time-varying systems with unobservable linearization, sufficient conditions are derived for the solvability of the observer design problem. Particularly, for the case where it is a priori known that the initial state of the system belongs to a given nonempty bounded subset of the state space, a Luenberger-type observer is constructed and for the general case we show that the state estimation is exhibited using a switching sequence of time-varying dynamics. The proposed methodology can be applied for further extensions; for instance, it can be shown that the observer design problem is solvable for the case $\dot{x}_i = f_i(t, x_1, ..., x_i) + a_i(t, x_1) r_i(t, x_{i+1})$, $i = 1, ..., n-1$, $\dot{x}_n = f_n(t, x_1, ..., x_n)$, $x_i \in \mathbb{R}$ with output $y := x_1$, under same hypotheses for the mappings $a_i$ and $f_i$ with those imposed in Theorem 1.1, together with forward completeness, and where $r_i(\cdot, \cdot)$, $i = 1, 2, ..., n-1$ above are $C^1$, in such a way that for each $t \geq 0$, the mappings $r_i(t, \cdot)$ are strictly monotone.